\documentclass{article}
\usepackage{amssymb}

\input{tcilatex}
\begin{document}

\begin{center}
\bigskip \textbf{Unique continuation properties for abstract Schrodinger
equations and applications}

\textbf{Veli Shakhmurov}

Department of Mechanical Engineering, Okan University, Akfirat, Tuzla 34959
Istanbul, Turkey,

E-mail: veli.sahmurov@okan.edu.tr

Institute of Mathematics and Mechanics, Azerbaijan National Academy of
Sciences, Azerbaijan, AZ1141, Baku, F. Agaev, 9,

E-mail: veli.sahmurov@gmail.com

A\textbf{bstract}
\end{center}

In this paper, Hardy's uncertainty principle and unique continuation
properties of Schrodinger equations with operator potentials in Hilbert
space-valued $L^{2}$ classes are obtained. Since the Hilbert space $H$ and
linear operators are arbitrary, by choosing the appropriate spaces and
operators we obtain numerous classes of Schrodinger type equations and its
finite and infinite many systems which occur in a wide variety of physical
systems.


\begin{center}
\bigskip\ \ \textbf{AMS 2010: 35Q41, 35K15, 47B25, 47Dxx, 46E40 }

\textbf{1. Introduction, definitions}
\end{center}

\bigskip In this paper, the unique continuation properties of the following
abstract Schr\"{o}dinger equation

\begin{equation}
i\partial _{t}u+\Delta u+A\left( x\right) u+V\left( x,t\right) u=0,\text{ }%
x\in R^{n},\text{ }t\in \left[ 0,T\right] ,  \tag{1.1}
\end{equation}%
are studied, where $A=A\left( x\right) $ is a linear and $V\left( x,t\right) 
$ is a given potential operator functions in a Hilbert space $H;$ $\Delta $
denotes the Laplace operator in $R^{n}$ and $u=$ $u(x,t)$ is the $H$-valued
unknown function. This linear result was then applied to show that two
regular solutions $u_{1}$, $u_{2}$ of non-linear Schr\"{o}dinger equations 
\begin{equation}
i\partial _{t}u+\Delta u+A\left( \left( x\right) \right) u=F\left( u,\bar{u}%
\right) ,\text{ }x\in R^{n},\text{ }t\in \left[ 0,T\right]  \tag{1.2}
\end{equation}
for general non-linearities $F$ must agree in $R^{n}\times \lbrack 0,T]$,
when $u_{1}-u_{2}$ and its gradient decay faster than any quadratic
exponential at times $0$ and $T$.

Hardy's uncertainty principle and unique continuation properties for Schr%
\"{o}dinger equations studied e.g in $\left[ \text{4-7}\right] $ and the
referances therein. Abstract differential equations studied e.g. in $\left[ 
\text{2, 9-19, 22, 24}\right] .$ However, there seems to be no such abstract
setting for nonlinear Schr\"{o}dinger equations except the local existence
of weak solution (cf. $\left[ \text{15}\right] .$ In contrast to these
results we will study the unique continuation properties of abstract Schr%
\"{o}dinger equations with the operator potentials.\ Since the Hilbert space 
$H$ is arbitrary and $A$ is a possible linear operator, by choosing $H$ and $%
A$ we can obtain numerous classes of Schr\"{o}dinger type equations and its
systems which occur in the different processes. Our main goal is to obtain
sufficient conditions on a solution $u$, the operator $A,$ potential $V$ and
the behavior of the solution at two different times $t_{0}$ and $t_{1}$
which guarantee that $u\left( x,t\right) \equiv 0$ for $x\in R^{n},$ $t\in %
\left[ 0,T\right] $. If we\ choose $H$ to be a concrete Hilbert space, for
example $H=L^{2}\left( \Omega \right) $, $A=L,$ where $\Omega $ is a domin
in $R^{m}$ with sufficientli smooth boundary and $L$ is an elliptic operator
then, we obtain the unique continuation properties of the followinng Schr%
\"{o}dinger equation%
\begin{equation}
\partial _{t}u=i\left( \Delta u+Lu\right) +V\left( x,t\right) u,\text{ }x\in
R^{n},\text{ }y\in \Omega ,\text{ }t\in \left[ 0,T\right] .  \tag{1.3}
\end{equation}%
\ Moreover, let we choose $H=L^{2}\left( 0,1\right) $ and $A$ to be
differential operator with Wentzell-Robin boundary condition defined by 
\begin{equation}
D\left( A\right) =\left\{ u\in W^{2,2}\left( 0,1\right) ,\text{ }Au\left(
j\right) =0,\text{ }j=0,1\right\} ,\text{ }  \tag{1.4}
\end{equation}%
\[
\text{ }A\left( x\right) u=a\left( x,y\right) u^{\left( 2\right) }+b\left(
x,y\right) u^{\left( 1\right) }, 
\]%
where $a,$ $b$ are sufficiently smooth functions on $R^{n}\times \left(
0,1\right) $ and $V\left( x,t\right) $ is a integral operator so that 
\[
V\left( x,t\right) u=\dint\limits_{0}^{1}K\left( x,\tau ,t\right) u\left(
x,y,\tau ,t\right) d\tau , 
\]%
where, $K=K\left( x,\tau ,t\right) $ is a complex valued bounded function.
From our general results we obtain the unique continuation properties of the
Wentzell-Robin type boundary value problem (BVP) for the following Schr\"{o}%
dinger equation 
\begin{equation}
\partial _{t}u=i\left( \Delta u+a\frac{\partial ^{2}u}{\partial y^{2}}+b%
\frac{\partial u}{\partial y}\right) +\dint\limits_{0}^{1}K\left( x,\tau
,t\right) u\left( x,y,\tau ,t\right) d\tau ,\text{ }  \tag{1.5}
\end{equation}%
\[
x\in R^{n},\text{ }y\in \left( 0,1\right) ,\text{ }t\in \left[ 0,T\right] , 
\]%
\ \ \ 

\begin{equation}
a\partial _{y}^{2}u\left( x,j,t\right) +b\partial _{y}u\left( x,j,t\right) =0%
\text{, }j=0,1.  \tag{1.6}
\end{equation}

Note that, the regularity properties of Wentzell-Robin type BVP for elliptic
equations were studied e.g. in $\left[ \text{11, 12 }\right] $ and the
references therein. Moreover, if put $H=l_{2}$ and choose $A$ to be a
infinite matrix $\left[ a_{mj}\right] $, $m,j=1,2,...,\infty ,$ then we
derive the unique continuation properties of the following system of Schr%
\"{o}dinger equation 
\begin{equation}
\partial _{t}u_{m}=i\left[ \Delta u_{m}+\sum\limits_{j=1}^{\infty }\left(
a_{mj}\left( x\right) +b_{mj}\left( x,t\right) \right) u_{j}\right] ,\text{ }%
x\in R^{n},\text{ }t\in \left( 0,T\right) ,  \tag{1.7}
\end{equation}%
where $a_{mj}$ are continuous and $b_{mj}$ are bounded functions.

Let $E$ be a Banach space. $L^{p}\left( \Omega ;E\right) $ denotes the space
of strongly measurable $E$-valued functions that are defined on the
measurable subset $\Omega \subset R^{n}$ with the norm

\[
\left\Vert f\right\Vert _{L^{p}}=\left\Vert f\right\Vert _{L^{p}\left(
\Omega ;E\right) }=\left( \int\limits_{\Omega }\left\Vert f\left( x\right)
\right\Vert _{E}^{p}dx\right) ^{\frac{1}{p}},\text{ }1\leq p<\infty \ . 
\]

Let $H$ be a Hilbert space and 
\[
\left\Vert u\right\Vert =\left\Vert u\right\Vert _{H}=\left( u,u\right)
_{H}^{\frac{1}{2}}=\left( u,u\right) ^{\frac{1}{2}}\text{ for }u\in H. 
\]

For $p=2$ and $E=H$, $L^{p}\left( \Omega ;E\right) $ becomes a $H$-valued
function space with inner product: 
\[
\left( f,g\right) _{L^{2}\left( \Omega ;H\right) }=\int\limits_{\Omega
}\left( f\left( x\right) ,g\left( x\right) \right) _{H}dx,\text{ }f,\text{ }%
g\in L^{2}\left( \Omega ;H\right) . 
\]

Here, $W^{s,2}\left( R^{n};H\right) $, $-\infty <s<\infty $ denotes the $H-$%
valued Sobolev space of order $s$ which is defined as: 
\[
W^{s,2}=W^{s,2}\left( R^{n};H\right) =\left( I-\Delta \right) ^{-\frac{s}{2}%
}L^{2}\left( R^{n};H\right) 
\]%
with the norm 
\[
\left\Vert u\right\Vert _{W^{s,2}}=\left\Vert \left( I-\Delta \right) ^{%
\frac{s}{2}}u\right\Vert _{L^{2}\left( R^{n};H\right) }<\infty . 
\]%
It clear that $W^{0,2}\left( R^{n};E\right) =L^{2}\left( R^{n};H\right) .$
Let $H_{0}$ and $H$ be two Hilbert spaces and $H_{0}$ is continuously and
densely embedded into $H$. Let $W^{s,2}\left( R^{n};H_{0},H\right) $ denote
the Sobolev-Lions type space, i.e., 
\[
W^{s,2}\left( R^{n};H_{0},H\right) =\left\{ u\in W^{s,2}\left(
R^{n};H\right) \cap L^{2}\left( R^{n};H_{0}\right) ,\right. \text{ } 
\]%
\[
\left. \left\Vert u\right\Vert _{W^{s,2}\left( R^{n};H_{0},H\right)
}=\left\Vert u\right\Vert _{L^{2}\left( R^{n};H_{0}\right) }+\left\Vert
u\right\Vert _{W^{s,2}\left( R^{n};H\right) }<\infty \right\} . 
\]

Let $C\left( \Omega ;E\right) $ denote the space of $E-$valued uniformly
bounded continious functions on $\Omega $ with norm 
\[
\left\Vert u\right\Vert _{C\left( \Omega ;E\right) }=\sup\limits_{x\in
\Omega }\left\Vert u\left( x\right) \right\Vert _{E}. 
\]

$C^{m}\left( \Omega ;E\right) $\ will denote the spaces of $E$-valued
uniformly bounded strongly continuous and $m$-times continuously
differentiable functions on $\Omega $ with norm 
\[
\left\Vert u\right\Vert _{C^{m}\left( \Omega ;E\right) }=\max\limits_{0\leq
\left\vert \alpha \right\vert \leq m}\sup\limits_{x\in \Omega }\left\Vert
D^{\alpha }u\left( x\right) \right\Vert _{E}. 
\]

Here, $O_{r}=\left\{ x\in R^{n},\text{ }\left\vert x\right\vert <r\right\} $
for $r>0$. Let $\mathbb{N}$ denote the set of all natural numbers, $\mathbb{C%
}$ denote the set of all complex numbers. Let $E_{1}$ and $E_{2}$ be two
Banach spaces. $B\left( E_{1},E_{2}\right) $ will denote the space of all
bounded linear operators from $E_{1}$ to $E_{2}.$ For $E_{1}=E_{2}=E$ it
will be denoted by $B\left( E\right) .$

\bigskip Here, $S=S(R^{n};E)$ denotes the $E-$valued Schwartz class, i.e.
the space of $E-$valued rapidly decreasing smooth functions on $R^{n},$
equipped with its usual topology generated by seminorms. $S(R^{n};\mathbb{C}%
) $ will be denoted by just $S$. Let $S^{\prime }(R^{n};E)$ denote the space
of all continuous linear operators, $L:S\rightarrow E$, equipped with
topology of bounded convergence.

Let $A=A\left( x\right) ,$ $x\in R^{n}$ be closed linear operator in $E$
with independent on $x\in R^{n}$ domain $D\left( A\right) $ that is dense on 
$E.$ The Fourier transformation of $A\left( x\right) ,$ i.e. $\hat{A}=FA=%
\hat{A}\left( \xi \right) $ is a linear operator defined as 
\[
\hat{A}\left( \xi \right) u\left( \varphi \right) =A\left( x\right) u\left( 
\hat{\varphi}\right) \text{ for }u\in S^{\prime }\left( R^{n};E\right) ,%
\text{ }\varphi \in S\left( R^{n}\right) . 
\]%
(For details see e.g $\left[ \text{1, Section 3}\right] $).

Let $\left[ A,B\right] $ be a commutator operator, i.e. 
\[
\left[ A,B\right] =AB-BA 
\]%
for linear operators $A$ and $B.$

Sometimes we use one and the same symbol $C$ without distinction in order to
denote positive constants which may differ from each other even in a single
context. When we want to specify the dependence of such a constant on a
parameter, say $\alpha $, we write $C_{\alpha }$.

\begin{center}
\textbf{2}. \textbf{Main results}
\end{center}

Let $A=A\left( x\right) ,$ $x\in R^{n}$ be closed linear operator in a
Hilbert space $H$ with independent on $x\in R^{n}$ domain $D\left( A\right) $
that is dense on $H.$ Let 
\[
H\left( A\right) =\left\{ u\in D\left( A\right) ,\text{ }\left\Vert
u\right\Vert _{H\left( A\right) }=\left\Vert Au\right\Vert _{H}+\left\Vert
u\right\Vert _{H}<\infty \right\} , 
\]
\[
X=L^{2}\left( R^{n};H\right) \text{, \ }X\left( A\right) =L^{2}\left(
R^{n};H\left( A\right) \right) ,\text{ }Y^{s}=W^{s,2}\left( R^{n};H\right) 
\text{, } 
\]%
\[
Y^{s}\left( A\right) =W^{s,2}\left( R^{n};H\left( A\right) \right) ,\text{ }%
B=L^{\infty }\left( R^{n};B\left( H\right) \right) \text{ and }\mu \left(
t\right) =\alpha t+\beta \left( 1-t\right) . 
\]

\textbf{Definition} \textbf{2.1}. A function $u\in L^{\infty }\left(
0,T;H\left( A\right) \right) $ is called a local weak solution to $\left(
1.1\right) $ on $\left( 0,T\right) $ if $u$ belongs to $L^{\infty }\left(
0,T;H\left( A\right) \right) \cap W^{1,2}\left( 0,T;H\right) $ and satisfies 
$\left( 1.1\right) $ in the sense of $L^{\infty }\left( 0,T;H\left( A\right)
\right) .$ In particular, if $\left( 0,T\right) $\ coincides with $\mathbb{R}
$, then $u$ is called a global weak solution to $\left( 1.1\right) .$

If the solution of $\left( 1.1\right) $ belongs to $C\left( \left[ 0,T\right]
;H\left( A\right) \right) \cap W^{2,2}\left( 0,T;H\right) ,$ then its called
a stronge solution.

Our main results in this paper is the following:

\textbf{Theorem 1. }Assume that the following condition are satisfied:

(1) $A=A\left( x\right) $ and $\frac{\partial A}{\partial x_{k}}$ are
symmetric operators in a Hilbert space $H$ with independent on $x\in R^{n}$
domain $D\left( \frac{\partial A}{\partial x_{k}}\right) =$ $D\left(
A\right) $ that is dense on $H.$ Moreover, 
\[
\dsum\limits_{k=1}^{n}\left( x_{k}\left[ A\frac{\partial f}{\partial x_{k}}-%
\frac{\partial A}{\partial x_{k}}f\right] ,f\right) _{X}\geq 0\text{ for }%
f\in L^{\infty }\left( 0,T;Y^{1}\left( A\right) \right) ; 
\]

(2) $A\left( x\right) A^{-1}\left( x_{0}\right) \in L^{1}\left(
R^{n};B\left( H\right) \right) $ for some $x_{0}\in R^{n}$\ and $V\left(
x,t\right) \in B\left( H\right) $ for $\left( x,t\right) \in R^{n}\times %
\left[ 0,1\right] .$ Moreover, there is a constant $C_{0}>0$ so that 
\[
\func{Im}\left( \left( A+V\right) \upsilon ,\upsilon \right) \geq C_{0}\eta
\left( x,t\right) \left\Vert \upsilon \left( x,t\right) \right\Vert ^{2} 
\]%
for $x\in R^{n},$ $t\in \left[ 0,T\right] $ and $\upsilon \in H$, where $%
\eta $ is a positive function in $L^{1}\left( 0,T;L^{\infty }\left(
R^{n}\right) \right) ;$

(3) either, $V\left( x,t\right) =V_{1}\left( x\right) +V_{2}\left(
x,t\right) $, where $V_{1}\left( x\right) \in B\left( H\right) $ for $x\in
R^{n}$ and%
\[
M_{1}=\sup\limits_{x\in R^{n}}\left\Vert V_{1}\left( x\right) \right\Vert
_{B\left( H\right) }<\infty ,\sup\limits_{t\in \left[ 0,1\right] }\left\Vert
e^{\left\vert x\right\vert ^{2}\mu ^{-2}\left( t\right) }V_{2}\left(
.,t\right) \right\Vert _{B}<\infty 
\]%
or%
\[
\lim\limits_{r\rightarrow \infty }\left\Vert V\right\Vert _{L^{1}\left(
0,1;L^{\infty }\left( R^{n}/O_{r}\right) ;B\left( H\right) \right) }=0. 
\]%
where 
\[
\text{ }\alpha ,\text{ }\beta >0,\text{ }\alpha \beta <2; 
\]

(4) $u\in C\left( \left[ 0,1\right] ;X\left( A\right) \right) $ is a
solution of the equation $\left( 1.1\right) $ and%
\[
\left\Vert e^{\beta ^{-2}\left\vert x\right\vert ^{2}}u\left( .,0\right)
\right\Vert _{X}<\infty ,\left\Vert e^{\alpha ^{-2}\left\vert x\right\vert
^{2}}u\left( .,1\right) \right\Vert _{X}<\infty . 
\]

Then $u\left( x,t\right) \equiv 0.$

As a direct consequence of Theorem 1 we get the following Hardy's
uncertainty principle result for the non-linear equations $(1.2)$.

\textbf{Theorem 2. }Assume that the assumptions (1)-(2) of Theorem 1 are
satisfied. Let $u_{1},u_{2}\in C\left( \left[ 0,1\right] ;Y^{k}\left(
A\right) \right) $, $k\in Z^{+}$ be stronge solutions of the equation $%
\left( 1.2\right) $ with $k>\frac{n}{2}.$ Moreover, assume $F\in C^{k}\left( 
\mathbb{C}^{2},\mathbb{C}\right) $ and $F\left( 0\right) =\partial
_{u}F\left( 0\right) =\partial _{\bar{u}}F\left( 0\right) =0.$ If there are $%
\alpha $, $\beta >0$ with $\alpha \beta <2$ such that 
\[
e^{-\beta ^{-2}\left\vert x\right\vert ^{2}}\left( u_{1}\left( .,0\right)
-u_{2}\left( .,0\right) \right) \in X\text{, }e^{-\alpha ^{-2}\left\vert
x\right\vert ^{2}}\left( u_{1}\left( .,1\right) -u_{2}\left( .,1\right)
\right) \in X 
\]%
then $u_{1}\equiv u_{2}.$

One of the results we get is the following one.

\textbf{Theorem 3. }Assume that the all conditions of Theorem 1 are
satisfied. Suppose $\Delta +A+V_{1}$ generates a bounded continious group.
Let $u\in C\left( \left[ 0,1\right] ;X\left( A\right) \right) $ be a
solution of $\left( 1.1\right) .$\ Then $\left\Vert e^{\left\vert
x\right\vert ^{2}\mu ^{-2}\left( t\right) }u\left( .,t\right) \right\Vert
_{X}^{\frac{1}{\mu \left( t\right) }}$ is logarithmically convex in $[0,1]$
and there is $N=N\left( \alpha ,\beta \right) $ such that 
\[
\left\Vert e^{\left\vert x\right\vert ^{2}\mu ^{-2}\left( t\right) }u\left(
.,t\right) \right\Vert _{X}^{\frac{1}{\mu \left( t\right) }}\leq 
\]%
\[
e^{N\left( M_{1}+M_{2}+M_{1}^{2}+M_{2}^{2}\right) }\left\Vert e^{\beta
^{-2}\left\vert x\right\vert ^{2}}u\left( .,0\right) \right\Vert _{X}^{\beta
\left( 1-t\right) \mu \left( t\right) }\left\Vert e^{\alpha ^{-2}\left\vert
x\right\vert ^{2}}u\left( .,1\right) \right\Vert _{X}^{\alpha t\mu \left(
t\right) }, 
\]%
when 
\[
M_{2}=e^{2B\left( V_{2}\right) }\sup\limits_{t\in \left[ 0,1\right]
}\left\Vert e^{\left\vert x\right\vert ^{2}\mu ^{-2}\left( t\right)
}V_{2}\left( .,t\right) \right\Vert _{B},\text{ }B\left( V_{2}\right)
=\sup\limits_{t\in \left[ 0,1\right] }\left\Vert \func{Re}V_{2}\left(
.,t\right) \right\Vert _{B}. 
\]%
Moreover,

\[
\sqrt{t\left( 1-t\right) }\left\Vert e^{\left\vert x\right\vert ^{2}\mu
^{-2}\left( t\right) }\nabla u\right\Vert _{L^{2}\left( R^{n}\times \left[
0,1\right] ;H\right) }\leq 
\]%
\[
e^{N\left( M_{1}+M_{2}+M_{1}^{2}+M_{2}^{2}\right) }\left[ \left\Vert
e^{\beta ^{-2}\left\vert x\right\vert ^{2}}u\left( .,0\right) \right\Vert
_{X}+\left\Vert e^{\alpha ^{-2}\left\vert x\right\vert ^{2}}u\left(
.,1\right) \right\Vert _{X}\right] . 
\]%
Here, we prove the following result for abstract parabolic equations with
variable coefficientes.

Consider the Cauchy problem for abstract parabolic equations with variable
operator coefficients 
\begin{equation}
\partial _{t}u=\Delta u+A\left( x\right) u+V\left( x,t\right) u,\text{ } 
\tag{2.1}
\end{equation}

\[
u\left( x,0\right) =f\left( x\right) ,\text{ }x\in R^{n},\text{ }t\in \left[
0,1\right] , 
\]%
where $A\left( x\right) $ is a linear and $V\left( x,t\right) $ is the given
potential operator functions in $H.$ By employing Theorem 1 we obtain

\textbf{Theorem 4. }Assume the assumptions (1)-(3) of Theorem 1 are
satisfied. Suppose $u\in L^{\infty }\left( 0,1;X\left( A\right) \right) \cap
L^{2}\left( 0,1;Y^{1}\right) $ is a solution of $\left( 2.1\right) $ and 
\[
\left\Vert f\right\Vert _{X}<\infty ,\left\Vert e^{\delta ^{-2}\left\vert
x\right\vert ^{2}}u\left( .,1\right) \right\Vert _{X}<\infty 
\]%
for some $\delta <1$. Then, $f\left( x\right) \equiv 0$ for $x\in R^{n}.$

First of all, we generalize the result G. H. Hardy (see e.g $\left[ 20\right]
$, p.131) about uncertainty principle for Fourier transform:

\textbf{Lemma 2.1. }Let \textbf{\ }$f\left( x\right) $ be $H$-valued
function for $x\in R^{n}$ and%
\[
\left\Vert f\left( x\right) \right\Vert =O\left( e^{-\frac{\left\vert
x\right\vert ^{2}}{\beta ^{2}}}\right) ,\text{ }\left\Vert \hat{f}\left( \xi
\right) \right\Vert =O\left( e^{-\frac{4\left\vert \xi \right\vert ^{2}}{%
\alpha ^{2}}}\right) ,\text{ }x\text{, }\xi \in R^{n}\text{ for }\alpha
\beta <4. 
\]

Then $f\left( x\right) \equiv 0.$ Also, if $\alpha \beta =4$ then$\left\Vert
f\left( x\right) \right\Vert $ is a constant multiple of $e^{-\frac{%
\left\vert x\right\vert ^{2}}{\beta ^{2}}}.$

\textbf{Proof. }Indeed, by employing Phragmen--Lindel\"{o}f theorem to the
classes of Hilbert-valued analytic functions and by reasoning as in $\left[ 8%
\right] $ we obtain the assertion.

Consider the Cauchy problem for free abstract Schr\"{o}dinger equation%
\begin{equation}
i\partial _{t}u+\Delta u+Au=0,\text{ }x\in R^{n},\text{ }t\in \left[ 0,1%
\right] ,  \tag{2.2}
\end{equation}

\[
u\left( x,0\right) =f\left( x\right) , 
\]%
where $A=A\left( x\right) $ is a linear operator in a Hilbert space $H$ with
independent on $x\in R^{n}$ domain $D\left( A\right) .$

The above result can be rewritten for solution of the $\left( 2.2\right) $
on $R^{n}\times \left( 0,\infty \right) $. Indeed, assume%
\[
\left\Vert u\left( x,0\right) \right\Vert =O\left( e^{-\frac{\left\vert
x\right\vert ^{2}}{\beta ^{2}}}\right) ,\text{ }\left\Vert u\left(
x,T\right) \right\Vert =O\left( e^{-\frac{\left\vert x\right\vert ^{2}}{%
\alpha ^{2}}}\right) \text{ for }\alpha \beta <4T. 
\]

Then $u\left( x,t\right) \equiv 0.$ Also, if $\alpha \beta =4T$, then $u$
has as a initial data a constant multiple of $e^{-\left( \frac{1}{\beta ^{2}}%
+\frac{i}{4T}\right) \left\vert y\right\vert ^{2}}.$

\textbf{Lemma 2.2. }Assume that $A$ is a symmetric operator in $H$ with
independent on $x\in R^{n}$ domain $D\left( A\right) $ that is dense on $H.$
Moreover, $A\left( x\right) A^{-1}\left( x_{0}\right) \in L^{1}\left(
R^{n};B\left( H\right) \right) $ for some $x_{0}\in R^{n}$. Then for $f\in
W^{s,2}\left( R^{n};H\right) ,$ $s\geq 0$ there is a generalized solution of 
$\left( 2.2\right) $ expressing as%
\begin{equation}
u\left( x,t\right) =F^{-1}\left[ e^{i\hat{A}_{\xi }t}\hat{f}\left( \xi
\right) \right] ,\text{ }\hat{A}_{\xi }=\hat{A}\left( \xi \right)
-\left\vert \xi \right\vert ^{2},  \tag{2.3}
\end{equation}

where $F^{-1}$ is the inverse Fourier transform and $\hat{A}\left( \xi
\right) $ denotes the Fourier transform of $A\left( x\right) .$

\textbf{Proof. }By applying the Fourier trasform to the problem $\left(
2.2\right) $ we get 
\begin{equation}
i\partial _{t}\hat{u}\left( \xi ,t\right) +\hat{A}_{\xi }\hat{u}\left( \xi
,t\right) =0,\text{ }x\in R^{n},\text{ }t\in \left[ 0,1\right] ,  \tag{2.4}
\end{equation}

\[
\hat{u}\left( \xi ,0\right) =\hat{f}\left( \xi \right) ,\text{ }\xi \in
R^{n}, 
\]

It is clear to see that the solution of the problem $\left( 2.4\right) $ can
be exspressed as 
\[
\hat{u}\left( \xi ,t\right) =e^{i\hat{A}_{\xi }t}\hat{f}\left( \xi \right) . 
\]

Hence, we obtain $\left( 2.3\right) .$

\begin{center}
\textbf{3. Estimates for solutions}
\end{center}

\bigskip We need the following lemmas for proving the main results. Consider
the abstract Schr\"{o}dinger equation%
\begin{equation}
\partial _{t}u=\left( a+ib\right) \left[ \Delta u+A\left( \left( x\right)
\right) u+V\left( x,t\right) u+F\left( x,t\right) \right] ,\text{ }x\in
R^{n},\text{ }t\in \left[ 0,1\right] ,  \tag{3.1}
\end{equation}%
where $a$, $b$ are real numbers, $A=A\left( x\right) $ is a linear operator, 
$V\left( x,t\right) $ is a given potential operator function in $H$ and $%
F\left( x,t\right) $ is a given $H$-valued function.

Let 
\[
\Phi \left( A,V\right) \upsilon =a\func{Re}\left( \left( A+V\right) \upsilon
,\upsilon \right) -b\func{Im}\left( \left( A+V\right) \upsilon ,\upsilon
\right) , 
\]%
\[
\text{ for }\upsilon =\upsilon \left( x,t\right) \in D\left( A\right) . 
\]

\textbf{Condition 3.1. }Assume that:

(1) $A=A\left( x\right) $ is a symmetric operator in Hilbert space $H$ with
independent on $x\in R^{n}$ domain $D\left( A\right) $ that is dense on $H;$

(2) $\frac{\partial A}{\partial x_{k}}$ are symmetric operators in $H$ with
independent on $x\in R^{n}$ domain $D\left( \frac{\partial A}{\partial x_{k}}%
\right) =$ $D\left( A\right) .$ Moreover, 
\[
\dsum\limits_{k=1}^{n}\left( x_{k}\left[ A\frac{\partial f}{\partial x_{k}}-%
\frac{\partial A}{\partial x_{k}}f\right] ,f\right) _{X}\geq 0,\text{ for }%
f\in L^{\infty }\left( 0,T;Y^{1}\left( A\right) \right) ; 
\]

\bigskip (3) there exists $d>0$ such that 
\begin{equation}
\left( A\upsilon ,\upsilon \right) \leq d\left\Vert \upsilon \left(
x,t\right) \right\Vert ^{2}\text{ for }\upsilon \in D\left( A\right) ; 
\tag{3.2}
\end{equation}

(4) $a>0,$ $b\in \mathbb{R};$ $V=V\left( x,t\right) \in B\left( H\right) $
and there is a constant $C_{0}>0$ so that 
\begin{equation}
\left\vert \Phi \left( A,V\right) \upsilon \left( x,t\right) \right\vert
\leq C_{0}\eta \left( x,t\right) \left\Vert \upsilon \left( .,t\right)
\right\Vert ^{2},  \tag{3.3}
\end{equation}%
for $x\in R^{n},$ $t\in \left[ 0,T\right] ,$ $T\in \left[ 0,1\right] $, $%
\upsilon \in D\left( A\right) $, where $\eta $ is a positive function in $%
L^{1}\left( 0,T;L^{\infty }\left( R^{n}\right) \right) $.

Let 
\[
\text{ }\left\vert \nabla \upsilon \right\vert
_{H}^{2}=\dsum\limits_{k=1}^{n}\left\Vert \frac{\partial \upsilon }{\partial
x_{k}}\right\Vert ^{2}\text{ for }\upsilon \in W^{1,2}\left( R^{n};H\right)
. 
\]

\textbf{Lemma 3.1. }Assume that the Condition 3.1 holds. Then the solution $%
u $ of $\left( 3.1\right) $ belonging to $L^{\infty }\left( 0,1;X\left(
A\right) \right) \cap L^{2}\left( 0,1;Y^{1}\right) $ satisfies the following
estimate%
\[
e^{M_{T}}\left\Vert e^{\phi \left( .,T\right) }u\left( .,T\right)
\right\Vert _{X}\leq M_{T}\left\Vert e^{\gamma \left\vert x\right\vert
^{2}}u\left( .,0\right) \right\Vert _{X}+\varkappa \left\Vert e^{\phi \left(
t\right) }F\right\Vert _{L^{1}\left( 0,T;X\right) }, 
\]%
where 
\[
\phi \left( x,t\right) =\frac{\gamma a\left\vert x\right\vert ^{2}}{%
a+4\gamma \left( a^{2}+b^{2}\right) t}\text{, }M_{T}=\left\Vert \eta
\right\Vert _{L^{1}\left( 0,T:L^{\infty }\left( R^{n}\right) \right) },\text{
}\varkappa =\sqrt{a^{2}+b^{2}},\text{ }\gamma \geq 0. 
\]

\textbf{Proof. }Let $\upsilon =e^{\varphi }u,$ where $\varphi $ is a
real-valued function to be chosen later. The function $\upsilon $ verifies%
\[
\partial _{t}\upsilon =S\upsilon +K\upsilon +\left( a+ib\right) \left[
\left( A+V\right) +e^{\varphi }F\right] ,\text{ }\left( x,t\right) \in
R^{n}\times \left[ 0,1\right] , 
\]
where $S$, $K$ are symmetric and skew-symmetric operators respectively given
by

\[
S=aA_{1}-ib\gamma B_{1}+\varphi _{t}+a\func{Re}V-b\func{Im}V,\text{ }%
K=ibA_{1}-a\gamma B_{1}+ 
\]%
\[
i\left( b\func{Re}\upsilon +a\func{Im}\upsilon \right) , 
\]%
here 
\[
A_{1}=\Delta +A\left( x\right) +\gamma ^{2}\left\vert \nabla \varphi
\right\vert ^{2},\text{ }B_{1}=2\nabla \varphi .\nabla +\Delta \varphi . 
\]%
By differentating inner product in $X$, we get 
\begin{equation}
\partial _{t}\left\Vert \upsilon \right\Vert _{X}^{2}=2\func{Re}\left(
S\upsilon ,\upsilon \right) _{X}+2\func{Re}\left( K\upsilon ,\upsilon
\right) _{X}+  \tag{3.4}
\end{equation}%
\[
2\func{Re}\left( \left( a+ib\right) e^{\upsilon }F,\upsilon \right) _{X}+2%
\func{Re}\left( \left( a+ib\right) \left( A+V\right) \upsilon ,\upsilon
\right) _{X},\text{ }t\geq 0. 
\]%
A formal integration by parts gives that%
\[
\func{Re}\left( S\upsilon ,\upsilon \right)
_{X}=-a\dint\limits_{R^{n}}\left\vert \nabla \upsilon \right\vert
_{H}^{2}dx+\dint\limits_{R^{n}}\left( a\gamma ^{2}\left\vert \nabla \varphi
\right\vert ^{2}+\varphi _{t}\right) \left\Vert \upsilon \right\Vert ^{2}dx+ 
\]%
\begin{equation}
a\dint\limits_{R^{n}}\left( A\upsilon ,\upsilon \right)
dx+a\dint\limits_{R^{n}}\left( \func{Re}V\upsilon ,\upsilon \right)
dx-b\dint\limits_{R^{n}}\left( \func{Im}V\upsilon ,\upsilon \right) dx, 
\tag{3.5}
\end{equation}%
\[
\func{Re}\left( K\upsilon ,\upsilon \right) _{X}=-a\gamma
\dint\limits_{R^{n}}\left[ \left( 2\nabla \varphi .\nabla \upsilon ,\upsilon
\right) +\Delta \varphi \left\Vert \upsilon \right\Vert ^{2}\right] dx, 
\]%
\[
\func{Re}\left( \left( a+ib\right) \left( A+V\right) \upsilon ,\upsilon
\right) _{X}=a\dint\limits_{R^{n}}\left( A\upsilon ,\upsilon \right) dx+a%
\func{Re}\dint\limits_{R^{n}}\left( V\upsilon ,\upsilon \right) dx- 
\]

\[
b\func{Im}\dint\limits_{R^{n}}\left( V\upsilon ,\upsilon \right)
dx=a\dint\limits_{R^{n}}\left( A\upsilon ,\upsilon \right)
dx+\dint\limits_{R^{n}}\Phi \left( A,V\right) \upsilon dx, 
\]

\[
\func{Re}\left( \left( a+ib\right) e^{\varphi }F,\upsilon \right) _{X}=a%
\func{Re}\dint\limits_{R^{n}}\left( e^{\varphi }F,\upsilon \right) dx-b\func{%
Im}\dint\limits_{R^{n}}\left( e^{\varphi }F,\upsilon \right) dx= 
\]%
\[
ae^{\varphi }\func{Re}\left( F,\upsilon \right) _{X}-be^{\varphi }\func{Im}%
\left( F,\upsilon \right) _{X}. 
\]

By using the Cauchy-Schwarz's inequality, by assumptions $\left( 3.1\right)
, $ $\left( 3.2\right) $, in view of $\left( 3.4\right) $ and $\left(
3.5\right) $ we obtain%
\[
\partial _{t}\left\Vert \upsilon \right\Vert ^{2}\leq 2\left\Vert a\func{Re}%
V-b\func{Im}V\right\Vert _{B}\left\Vert \upsilon \right\Vert
_{X}^{2}+2\varkappa \left\Vert e^{\varphi }F\left( t,.\right) \right\Vert
_{X}\left\Vert \upsilon \right\Vert _{X}+d\left\Vert \upsilon \right\Vert
_{X}^{2}, 
\]%
where $a,$ $b$ and $\varphi $ are such that%
\begin{equation}
\left( \text{ }a+\frac{b^{2}}{a}\right) \left\vert \nabla \varphi
\right\vert ^{2}+\varphi _{t}\leq 0\text{ in }R_{+}^{n+1}.  \tag{3.6}
\end{equation}%
The remainig part of the proof is obtained by reasoning as in $\left[ \text{%
7, Lemma 1}\right] .$

When $\varphi \left( x,t\right) =q\left( t\right) \psi \left( x\right) $, it
suffices that%
\begin{equation}
\left( a+\frac{b^{2}}{a}\right) q^{2}\left( t\right) \left\vert \nabla \psi
\right\vert ^{2}+q^{\prime }\left( t\right) \psi \left( x\right) \leq 0. 
\tag{3.7}
\end{equation}%
If we put $\psi \left( x\right) =\left\vert x\right\vert ^{2}$ then $\left(
3.7\right) $ holds, when 
\begin{equation}
q^{\prime }\left( t\right) =-4\left( a+\frac{b^{2}}{a}\right) q^{2}\left(
t\right) ,\text{ }q\left( 0\right) =\gamma ,\text{ }\gamma \geq 0.  \tag{3.8}
\end{equation}

Let 
\[
\psi _{r}\left( x\right) =\left\{ 
\begin{array}{c}
\left\vert x\right\vert ^{2}\text{, }\left\vert x\right\vert <r \\ 
\infty \text{, \ \ \ \ }\left\vert x\right\vert >r\text{\ \ \ \ \ \ \ \ }%
\end{array}%
.\right. 
\]%
Regularize $\psi _{r}$ with a radial mollifier $\theta _{\rho }$ and set 
\[
\varphi _{\rho ,r}\left( x,t\right) =q\left( t\right) \theta _{\rho }\ast
\psi _{r}\left( x\right) ,\text{ }\upsilon _{\rho ,r}\left( x,t\right)
=e^{\varphi _{\rho ,r}}u,\text{ } 
\]%
where $q\left( t\right) =\gamma a\left[ a+4\gamma \left( a^{2}+b^{2}\right) t%
\right] ^{-1}$ is the solution to $\left( 3.5\right) $. Because the right
hand side of $(3.5)$ only involves the first derivatives of $\varphi $, $%
\psi _{r}$ is Lipschitz and bounded at infinity, 
\[
\theta _{\rho }\ast \psi _{r}\left( x\right) \leq \theta _{\rho }\ast
\left\vert x\right\vert ^{2}=C\left( n\right) \rho ^{2} 
\]%
and $(3.6)$ holds uniformly in $\rho $ and $r$, when $\varphi $ is replaced
by $\varphi _{\rho ,r}$. Hence, it follows that the estimate%
\[
e^{M_{T}}\left\Vert e^{\phi \left( T\right) }u\left( T\right) \right\Vert
_{X}\leq M_{T}\left\Vert e^{\gamma \left\vert x\right\vert ^{2}}u\left(
0\right) \right\Vert _{X}+\sqrt{a^{2}+b^{2}}\left\Vert e^{\varphi _{\rho
,r}}F\right\Vert _{L^{1}\left( 0,T;X\right) } 
\]%
holds uniformly in $\rho $ and $r.$ The assertion is obtained after letting $%
\rho $ tend to zero and $r$ to infinity.

\textbf{Remark 3.1. }It should be noted that if $H=\mathbb{C}$, $A=0$ and $%
V\left( x,t\right) $ is a complex valued function, then the abstract
condition $\left( 3.3\right) $ can be replaced by 
\[
M_{T}=\left\Vert a\func{Re}V-b\func{Im}V\right\Vert _{L^{1}\left(
0,T;L^{\infty }\left( R^{n}\right) \right) }<\infty . 
\]%
Moreover, if $A\left( x\right) $ and $V\left( x,t\right) $ for $x\in R^{n},$ 
$t\in \left[ 0,T\right] $ are bounded operators in $H$, then by using
Cauchy-Schwarz's inequality, the assumption $\left( 3.3\right) $ becames as: 
\[
\left\vert \Phi \left( A,V\right) \upsilon \right\vert \leq \varkappa
\left\Vert \left( A+V\right) \upsilon \right\Vert _{H}\left\Vert \upsilon
\right\Vert _{H}\leq \varkappa \left\Vert A+V\right\Vert _{B\left( H\right)
}\left\Vert \upsilon \right\Vert ^{2}. 
\]

Let 
\[
\text{ }Q\left( t\right) =\left( f,f\right) _{X}\text{, }D\left( t\right)
=\left( Sf,f\right) _{X},\text{ }N\left( t\right) =D\left( t\right)
Q^{-1}\left( t\right) ,\text{ }\partial _{t}S=S_{t}. 
\]

\textbf{Lemma 3.2. }Assume\ that $S=S\left( t\right) $ is a symmetric, $%
K=K\left( t\right) $ is a skew-symmetric operators in $H$, $G\left(
x,t\right) $ is a positive funtion and $f(x,t)$ is a reasonable function.
Then, 
\[
Q^{^{\prime \prime }}\left( t\right) =2\partial _{t}\func{Re}\left( \partial
_{t}f-Sf-Kf,f\right) _{X}+2\left( S_{t}f+\left[ S,K\right] f,f\right) _{X}+ 
\]

\begin{equation}
\left\Vert \partial _{t}f-Sf+Kf\right\Vert _{X}^{2}-\left\Vert \partial
_{t}f-Sf-Kf\right\Vert _{X}^{2}  \tag{3.9}
\end{equation}%
and 
\[
\partial _{t}N\left( t\right) \geq Q^{-1}\left( t\right) \left[ \left(
S_{t}f+\left[ S,K\right] f,f\right) _{X}-\frac{1}{2}\left\Vert \partial
_{t}f-Sf-Kf\right\Vert _{X}^{2}\right] . 
\]

\bigskip Moreover, if 
\[
\left\Vert \partial _{t}f-Sf-Kf\right\Vert _{H}\leq M_{1}\left\Vert
f\right\Vert _{H}+G\left( x,t\right) ,\text{ }S_{t}+\left[ S,K\right] \geq
-M_{0}\text{ } 
\]%
for $x\in R^{n}$, $t\in \left[ 0,1\right] $ and%
\[
M_{2}=\sup\limits_{t\in \left[ 0,1\right] }\left\Vert G\left( .,t\right)
\right\Vert _{L^{2}\left( R^{n}\right) }\left( \left\Vert f\left( .,t\right)
\right\Vert _{X}\right) ^{-1}<\infty . 
\]

Then $Q\left( t\right) $ is logarithmically convex in $[0,1]$ and there is a
constant $M$ such that 
\[
Q\left( t\right) \leq e^{M\left(
M_{0}+M_{1}+M_{2}+M_{1}^{2}+M_{2}^{2}\right) }Q^{1-t}\left( 0\right)
Q^{t}\left( 1\right) \text{, }0\leq t\leq 1. 
\]

\textbf{Proof.} The lemma is verifying in a similar way as in $\left[ \text{%
7, Lemma 2}\right] $ by replacing the inner product and norm of $L^{2}\left(
R^{n}\right) $ with inner product and norm of the space $L^{2}\left(
R^{n};H\right) .$

\bigskip \textbf{Lemma 3.3. }Assume that the Condition 3.1 holds. Moreover,
suppose%
\[
\sup\limits_{t\in \left[ 0,1\right] }\left\Vert V\left( .,t\right)
\right\Vert _{B}\leq M_{1}\text{, }\left\Vert e^{\gamma \left\vert
x\right\vert ^{2}}u\left( .,0\right) \right\Vert _{X}<\infty ,\text{ } 
\]%
\[
\left\Vert e^{\gamma \left\vert x\right\vert ^{2}}u\left( .,1\right)
\right\Vert _{X}<\infty ,\text{ }M_{2}=\sup\limits_{t\in \left[ 0,1\right] }%
\frac{\left\Vert e^{\gamma \left\vert x\right\vert ^{2}}F\left( .,t\right)
\right\Vert _{X}}{\left\Vert u\right\Vert _{X}}<\infty . 
\]

Then, for solution $u\in L^{\infty }\left( 0,1;X\left( A\right) \right) \cap
L^{2}\left( 0,1;Y^{1}\right) $ of $\left( 3.1\right) $, $e^{\gamma
\left\vert x\right\vert ^{2}}u\left( .,t\right) $ is logarithmically convex
in $[0,1]$ and there is a constant $N$ such that%
\begin{equation}
\left\Vert e^{\gamma \left\vert x\right\vert ^{2}}u\left( .,t\right)
\right\Vert _{X}\leq e^{NM\left( a,b\right) }\left\Vert e^{\gamma \left\vert
x\right\vert ^{2}}u\left( .,0\right) \right\Vert _{X}^{1-t}\left\Vert
e^{\gamma \left\vert x\right\vert ^{2}}u\left( .,1\right) \right\Vert
_{X}^{t}  \tag{3.10}
\end{equation}%
where 
\[
M\left( a,b\right) =\varkappa ^{2}\left( \gamma M_{1}^{2}+M_{2}^{2}\right)
+\varkappa \left( M_{1}+M_{2}\right) 
\]%
when $0\leq t\leq 1.$

\textbf{Proof. }Let $f=e^{\gamma \varphi }u,$ where $\varphi $ is a
real-valued function to be chosen. The function $f\left( x\right) $ verifies%
\begin{equation}
\partial _{t}f=Sf+Kf+\left( a+ib\right) \left( Vf+e^{\gamma \varphi
}F\right) \text{ in }R^{n}\times \left[ 0,1\right] \text{,}  \tag{3.11}
\end{equation}%
where $S$, $K$ are symmetric and skew-symmetric operator, respectively given
by 
\[
S=aA_{1}-ib\gamma B_{1}+\varphi _{t}+a\func{Re}V-b\func{Im}V,\text{ }%
K=ibA_{1}-a\gamma B_{1}+ 
\]%
\begin{equation}
i\left( b\func{Re}\upsilon +a\func{Im}\upsilon \right) ,  \tag{3.12}
\end{equation}%
here 
\[
A_{1}=\Delta +A\left( x\right) +\gamma ^{2}\left\vert \nabla \varphi
\right\vert ^{2},\text{ }B_{1}=2\nabla \varphi .\nabla +\Delta \varphi . 
\]

\bigskip A calculation shows that, 
\begin{equation}
S_{t}+\left[ S,K\right] =\gamma \partial _{t}^{2}\varphi +2\gamma
^{2}a\nabla \varphi .\nabla \varphi _{t}-2ib\gamma \left( 2\nabla \varphi
_{t}.\nabla +\Delta \varphi _{t}\right) -  \tag{3.13}
\end{equation}%
\[
\gamma \varkappa ^{2}\left[ 4\nabla .\left( D^{2}\varphi \nabla \right)
-4\gamma ^{2}D^{2}\varphi \nabla \varphi +\Delta ^{2}\varphi \right] +2\left[
A\left( x\right) \nabla \varphi .\nabla -\nabla \varphi .\nabla A\right] . 
\]%
If we put $\varphi =\left\vert x\right\vert ^{2}$, then $\left( 3.13\right) $
reduce the following%
\[
S_{t}+\left[ S,K\right] =-\gamma \varkappa ^{2}\left[ 8\Delta -32\gamma
^{2}\left\vert x\right\vert ^{2}\right] +\frac{d}{dt}A+2\left[ A\left(
x\right) \nabla \varphi .\nabla -\nabla \varphi .\nabla A\right] . 
\]

\bigskip Moreover by assumtion (2), 
\[
\left( S_{t}f+\left[ S,K\right] f,f\right) =\gamma \varkappa
\dint\limits_{R^{n}}\left( 8\left\vert \nabla f\right\vert _{H}^{2}+32\gamma
^{2}\left\vert x\right\vert ^{2}\left\Vert f\right\Vert ^{2}\right) dx+ 
\]%
\begin{equation}
2\dint\limits_{R^{n}}\left( \left[ A\left( x\right) \nabla \varphi .\nabla
f-\nabla \varphi .\nabla Af\right] ,f\right) dx\geq 0.  \tag{3.14}
\end{equation}

This identity, the condition on $V$ and $\left( 3.14\right) $ imply that%
\begin{equation}
\left\Vert \partial _{t}f-Sf-Kf\right\Vert _{X}\leq \varkappa ^{2}\left(
M_{1}\left\Vert f\right\Vert _{X}+e^{\gamma \varphi }\left\Vert F\right\Vert
_{X}\right) \text{.}  \tag{3.15}
\end{equation}%
If we knew that the quantities and calculations involved in the proof of
Lemma 3.2 (similar as in $\left[ \text{7, Lemma 2}\right] $) were finite and
correct, when $f=e^{\gamma \left\vert x\right\vert ^{2}}u$ we would have the
logarithmic convexity of $Q\left( t\right) =\left\Vert e^{\gamma \left\vert
x\right\vert ^{2}}u\left( .,t\right) \right\Vert _{X}$ \ and the estimate $%
\left( 3.10\right) $ from Lemma 3.2. But this fact is verifying by reasonong
as in $\left[ \text{7, Lemma 3}\right] .$

Let 
\[
\sigma =\sqrt{t\left( 1-t\right) }e^{\gamma \left\vert x\right\vert ^{2}},%
\text{ }Z=L^{2}\left( \left[ 0,1\right] \times R^{n};H\right) . 
\]

\textbf{Lemma 3.4. }Assume that $a$, $b$, $u$, $A$ and $V$ are as in Lemma
3.3 and $\gamma >0$. Then,

\ \ \ \ \ \ \ \ \ \ \ \ \ \ \ \ \ \ \ \ \ \ \ \ \ \ \ \ \ \ \ \ \ \ \ \ \ \
\ \ \ \ \ \ \ \ \ \ \ \ \ \ \ \ \ \ \ \ \ \ \ \ \ \ \ \ \ \ \ \ \ \ \ \ \ \
\ \ \ \ \ \ \ \ \ \ \ \ \ \ \ 
\[
\left\Vert \sigma \nabla u\right\Vert _{Z}+\left\Vert \sigma \left\vert
x\right\vert u\right\Vert _{Z}\leq N\left[ \left( 1+M_{1}\right) \right] %
\left[ \sup\limits_{t\in \left[ 0,1\right] }\left\Vert e^{\gamma \left\vert
x\right\vert ^{2}}u\left( .,t\right) \right\Vert _{X}+\sup\limits_{t\in %
\left[ 0,1\right] }\left\Vert e^{\gamma \left\vert x\right\vert ^{2}}F\left(
.,t\right) \right\Vert _{Z}\right] \text{,} 
\]%
where $N$ is bounded number, when $\gamma \ $and $\varkappa $ are bounded
below.

\textbf{Proof. }The integration by parts shows that 
\[
\dint\limits_{R^{n}}\left( \left\vert \nabla f\right\vert _{H}^{2}+4\gamma
^{2}\left\vert x\right\vert ^{2}\left\Vert f\right\Vert \right)
dx=\dint\limits_{R^{n}}\left[ e^{2\gamma \left\vert x\right\vert ^{2}}\left(
\left\vert \nabla u\right\vert _{H}^{2}-2n\gamma \right) \left\Vert
u\right\Vert ^{2}\right] dx, 
\]%
when $f=e^{\gamma \left\vert x\right\vert ^{2}}u$, while integration by
parts, the Cauchy-Schwarz's inequality and the identity, $n=\nabla $%
\textperiodcentered $x$, give that 
\[
\dint\limits_{R^{n}}\left( \left\vert \nabla f\right\vert _{H}^{2}+4\gamma
^{2}\left\vert x\right\vert ^{2}\left\Vert f\right\Vert \right) dx\geq
2\gamma n\left\Vert f\right\Vert _{X}^{2}\text{.} 
\]

The sum of the last two formulae gives the inequality%
\begin{equation}
2\dint\limits_{R^{n}}\left( \left\vert \nabla f\right\vert _{H}^{2}+4\gamma
^{2}\left\vert x\right\vert ^{2}\left\Vert f\right\Vert \right) dx\geq
\dint\limits_{R^{n}}e^{\gamma \left\vert x\right\vert ^{2}}\left\vert \nabla
f\right\vert _{H}^{2}dx\text{.}  \tag{3.16}
\end{equation}%
Integration over $[0,1]$ of $t(1-t)$ times the formula (3.6) for $%
Q^{^{\prime \prime }}\left( t\right) $ and integration by parts, shows that 
\begin{equation}
2\dint\limits_{0}^{1}t(1-t)\left( S_{t}f+\left[ S,K\right] f,f\right)
_{X}dt+\dint\limits_{0}^{1}Q\left( t\right) dt\leq Q\left( 1\right) +Q\left(
0\right) +  \tag{3.17}
\end{equation}

\[
2\dint\limits_{0}^{1}1-2t)\func{Re}\left( \partial _{t}f-Sf-Kf,f\right)
_{X}dx+\dint\limits_{0}^{1}t(1-t)\left\Vert \partial _{t}f-Sf-Kf\right\Vert
_{X}^{2}dt. 
\]

Assuming again that the last two calculations are justified for $f=e^{\gamma
\left\vert x\right\vert ^{2}}.$ Then $\left( 3.14\right) -\left( 3.17\right) 
$ implay the assertion.

\begin{center}
\textbf{4. Appell transformation in abstract functon spaces}
\end{center}

\bigskip Let 
\[
\rho \left( t\right) =\alpha \left( 1-t\right) +\beta t,\text{ }\varphi
\left( x,t\right) =\frac{\left( \alpha -\beta \right) \left\vert
x^{2}\right\vert }{4\left( a+ib\right) \rho \left( t\right) },\text{ } 
\]%
\[
\nu \left( s\right) =\left[ \gamma \alpha \beta \rho ^{2}\left( s\right) +%
\frac{\left( \alpha -\beta \right) a}{4\left( a^{2}+b^{2}\right) }\rho
\left( s\right) \right] . 
\]

\bigskip \textbf{Lemma 4.1. }Assume $A$ and $V$ are as in Lemma 3.3 and $%
u=u\left( x,s\right) $ is a solution of the equation 
\[
\partial _{s}u=\left( a+ib\right) \left[ \Delta u+Au+V\left( y,s\right)
u+F\left( y,s\right) \right] ,\text{ }y\in R^{n},\text{ }s\in \left[ 0,1%
\right] . 
\]%
Let $a+ib\neq 0$, $\gamma \in \mathbb{R}$ and $\alpha $, $\beta \in \mathbb{R%
}_{+}$. Set 
\begin{equation}
\tilde{u}\left( x,t\right) =\left( \sqrt{\alpha \beta }\rho ^{-1}\left(
t\right) \right) ^{\frac{n}{2}}u\left( \sqrt{\alpha \beta }x\rho ^{-1}\left(
t\right) ,\beta t\rho ^{-1}\left( t\right) \right) e^{\varphi }.  \tag{4.1}
\end{equation}

\bigskip Then, $\tilde{u}\left( x,t\right) $\ verifies the equation 
\[
\partial _{t}\tilde{u}=\left( a+ib\right) \left[ \Delta \tilde{u}+A\tilde{u}+%
\tilde{V}\left( x,t\right) u+\tilde{F}\left( x,t\right) \right] ,\text{ }%
x\in R^{n},\text{ }t\in \left[ 0,1\right] 
\]%
with 
\[
\tilde{V}\left( x,t\right) =\alpha \beta \rho ^{-2}\left( t\right) V\left( 
\sqrt{\alpha \beta }x\rho ^{-1}\left( t\right) ,\beta t\rho ^{-1}\left(
t\right) \right) , 
\]

\[
\text{ }\tilde{F}\left( x,t\right) =\left( \sqrt{\alpha \beta }\rho
^{-1}\left( t\right) \right) ^{\frac{n}{2}+2}\left( \sqrt{\alpha \beta }%
x\rho ^{-1}\left( t\right) ,\beta t\rho ^{-1}\left( t\right) \right) . 
\]%
Moreover, 
\[
\left\Vert e^{\gamma \left\vert x\right\vert ^{2}}\tilde{F}\left( .,t\right)
\right\Vert _{X}=\alpha \beta \rho ^{-2}\left( t\right) e^{\nu \left\vert
y\right\vert ^{2}}\left\Vert F\left( s\right) \right\Vert _{X}\text{ and }%
\left\Vert e^{\gamma \left\vert x\right\vert ^{2}}\tilde{u}\left( .,t\right)
\right\Vert _{X}=e^{\nu \left\vert y\right\vert ^{2}}\left\Vert u\left(
s\right) \right\Vert _{X} 
\]%
when $s=\mu \left( t\right) $ and $\gamma \in \mathbb{R}$.

\textbf{Proof. }If $u$ is a solution of the equation%
\begin{equation}
\partial _{s}u=\left( a+ib\right) \left[ \Delta u+Au+Q\left( y,s\right) %
\right] ,\text{ }y\in R^{n},\text{ }s\in \left[ 0,1\right]  \tag{4.2}
\end{equation}%
then, the function $u_{1}\left( x,t\right) =u\left( \sqrt{r}x,rt+\tau
\right) $ verifies%
\[
\partial _{t}u_{1}=\left( a+ib\right) \left[ \Delta u_{1}+Au_{1}+rQ\left( 
\sqrt{r}x,rt+\tau \right) \right] ,\text{ }y\in R^{n},\text{ }s\in \left[ 0,1%
\right] 
\]%
and $u_{2}\left( x,t\right) =t^{-\frac{n}{2}}u\left( \frac{x}{t},\frac{1}{t}%
\right) e^{\frac{\left\vert x\right\vert ^{2}}{4\left( a+ib\right) t}}$ is a
\ solution to%
\[
\partial _{t}u_{2}=-\left( a+ib\right) \left[ \Delta u_{2}+Au+t^{-\left( 2+%
\frac{n}{2}\right) }Q\left( \frac{x}{t},\frac{1}{t}\right) e^{\frac{%
\left\vert x\right\vert ^{2}}{4\left( a+ib\right) t}}\right] ,\text{ }y\in
R^{n},\text{ }s\in \left[ 0,1\right] . 
\]%
These two facts and the sequel of changes of variables below verifies the
Lemma, when $\alpha >\beta ,$ i.e.%
\[
u\left( \sqrt{\frac{\alpha \beta }{\alpha -\beta }}x,\frac{\alpha \beta }{%
\alpha -\beta }t-\frac{\beta }{\alpha -\beta }\right) 
\]%
\ \ is a solution to the same non-homogeneous equation but with right-hand
side 
\[
\frac{\alpha \beta }{\alpha -\beta }Q\left( \sqrt{\frac{\alpha \beta }{%
\alpha -\beta }}x,\frac{\alpha \beta }{\alpha -\beta }t-\frac{\beta }{\alpha
-\beta }\right) . 
\]

The function, 
\[
\frac{1}{\left( \alpha -t\right) ^{\frac{n}{2}}}u\left( \frac{\sqrt{\alpha
\beta }x}{\sqrt{\alpha -\beta }\left( \alpha -t\right) },\frac{\alpha \beta 
}{\left( \alpha -\beta \right) \left( \alpha -t\right) }-\frac{\beta }{%
\alpha -\beta }\right) e^{\frac{\left\vert x\right\vert ^{2}}{4\left(
a+ib\right) \left( \alpha -t\right) }} 
\]%
verifies $(4.2)$ with right-hand side%
\[
\frac{\alpha \beta }{\left( \alpha -\beta \right) \left( \alpha -t\right) ^{%
\frac{n}{2}+2}}Q\left( \frac{\sqrt{\alpha \beta }x}{\sqrt{\alpha -\beta }%
\left( \alpha -t\right) },\frac{\alpha \beta }{\left( \alpha -\beta \right)
\left( \alpha -t\right) }-\frac{\beta }{\alpha -\beta }\right) e^{\frac{%
\left\vert x\right\vert ^{2}}{4\left( a+ib\right) \left( \alpha -t\right) }%
}. 
\]%
Replacing $\left( x,t\right) $ by $\left( \sqrt{\alpha -\beta }x,\left(
\alpha -\beta \right) t\right) $ we get that%
\begin{equation}
\rho ^{-\frac{n}{2}}\left( t\right) u\left( \sqrt{\alpha \beta }\rho
^{-1}\left( t\right) x,\frac{\alpha \beta \rho ^{-1}\left( t\right) }{\left(
\alpha -\beta \right) }-\frac{\beta }{\alpha -\beta }\right) e^{^{\left(
\alpha -\beta \right) }\frac{\left\vert x\right\vert ^{2}\rho \left(
t\right) }{4\left( a+ib\right) }}  \tag{4.3}
\end{equation}%
is a solution of $\left( 4.2\right) $ but with right-hand%
\begin{equation}
\rho ^{-\left( \frac{n}{2}+2\right) }\left( t\right) Q\left( \sqrt{\alpha
\beta }\rho ^{-1}\left( t\right) x,\frac{\alpha \beta \rho ^{-1}\left(
t\right) }{\left( \alpha -\beta \right) }-\frac{\beta }{\alpha -\beta }%
\right) e^{^{\left( \alpha -\beta \right) }\frac{\left\vert x\right\vert
^{2}\rho \left( t\right) }{4\left( a+ib\right) }}.  \tag{4.4}
\end{equation}%
Finally, observe that%
\[
s=\beta t\rho \left( t\right) =\frac{\alpha \beta \rho ^{-1}\left( t\right) 
}{\left( \alpha -\beta \right) }-\frac{\beta }{\alpha -\beta } 
\]%
and multiply $\left( 4.3\right) $ and $\left( 4.4\right) $ we obtain the
assertion for $\alpha >\beta .$ The case $\beta >\alpha $ follows by
reversing by changes of variables, $s^{\prime }=1-s$ and $t^{\prime }=1-t.$

\begin{center}
\bigskip \textbf{\ 5. Variable coefficients. Proof of Theorem 3}
\end{center}

We are ready to prove Theorem 3.\textbf{\ }Let 
\[
B=L^{1}\left( 0,1;L^{\infty }\left( R^{n};B\left( H\right) \right) \right) ,%
\text{ }B\left( R\right) =L^{1}\left( 0,1;L^{\infty }\left(
R^{n}/O_{r};B\left( H\right) \right) \right) . 
\]

\textbf{Proof of Theorem 3. }We may assume that $\alpha \neq \beta $. The
case $\alpha =\beta $ follows from the latter by replacing $\beta $ by $%
\beta +\delta ,$ $\delta >0$, and letting $\delta $ tend to zero. We may
also assume that $\alpha <\beta $. Otherwise, replace $u$ by $\bar{u}(1-t)$.
Assume $a>0.$ Set $W=\Delta +A+V_{1}.$ By Lemma 2.2 the problem

\begin{equation}
\partial _{t}u=\left( a+ib\right) \left[ \Delta u+A\left( x\right)
u+V_{1}\left( x\right) u\right] ,\text{ }x\in R^{n},\text{ }t\in \left[ 0,1%
\right] ,  \tag{5.1}
\end{equation}%
\[
u\left( x,0\right) =u_{0}\left( x\right) . 
\]%
has a solution $u=U\left( t\right) u_{0}=e^{t\left( a+ib\right) W}u_{0}\in
C\left( \left[ 0,1\right] ;X\left( A\right) \right) $, \ where%
\[
U\left( t\right) =F^{-1}\left[ e^{-Q\left( \xi \right) }\right] ,\text{ }%
Q\left( \xi \right) =\left( a+ib\right) \left[ -\left\vert \xi \right\vert
^{2}+\hat{A}\left( \xi \right) +\hat{V}_{1}\left( \xi \right) \right] , 
\]%
here, $F^{-1}$ is the inverse Fourier transform, $\hat{A}\left( \xi \right)
, $ $\hat{V}_{1}\left( \xi \right) $ respectively denote the Fourier
transforms of $A\left( x\right) ,$ $V_{1}\left( x\right) .$ By reasoning as
\ the Duhamel principle we get that the problem 
\[
\partial _{t}u=\left( a+ib\right) \left[ \Delta u+A\left( x\right) u+V\left(
x,t\right) \right] u,\text{ }x\in R^{n},\text{ }t\in \left[ 0,1\right] , 
\]%
\[
u\left( x,0\right) =u_{0}\left( x\right) . 
\]%
has a solution expressing as 
\begin{equation}
u\left( x,t\right) =H\left( t\right) u_{0}+i\dint\limits_{0}^{t}H\left(
t-s\right) V_{2}\left( x,s\right) u\left( x,s\right) ds\text{ }  \tag{5.2}
\end{equation}%
\[
\text{for }x\in R^{n},\text{ }s\in \left[ 0,1\right] , 
\]%
where 
\[
e^{itW}=H\left( t\right) =H\left( t,x\right) =F^{-1}\left[ e^{iQ\left( \xi
\right) }\right] . 
\]

For $0\leq \varepsilon \leq 1$ set 
\begin{equation}
F_{\varepsilon }\left( x,t\right) =\frac{i}{\varepsilon +i}e^{\varepsilon
tW}V_{2}\left( x,t\right) u\left( x,t\right)  \tag{5.3}
\end{equation}%
and 
\begin{equation}
u_{\varepsilon }\left( x,t\right) =e^{\left( \varepsilon +i\right)
tW}u_{0}+\left( \varepsilon +i\right) \dint\limits_{0}^{t}e^{\left(
\varepsilon +i\right) \left( t-s\right) W}F_{\varepsilon }\left( x,s\right)
u\left( x,s\right) ds.\text{ }  \tag{5.4}
\end{equation}

Then, $u_{\varepsilon }\left( x,t\right) \in L^{\infty }\left( 0,1;X\left(
A\right) \right) \cap L^{2}\left( R^{n};Y^{1}\right) $ and satisfies 
\[
\partial _{t}u_{\varepsilon }=\left( \varepsilon +i\right) \left(
Wu+F_{\varepsilon }\right) \text{ in }R^{n}\times \left[ 0,1\right] , 
\]%
\[
u_{\varepsilon }\left( .,0\right) =u_{0}\left( .\right) . 
\]

The identities 
\[
e^{\left( z_{1}+z_{2}\right) W}=e^{z_{1}W}e^{z_{2}W}\text{, when }\func{Re}%
z_{1}\text{, }\func{Re}z_{2}\geq 0, 
\]%
and $\left( 5.2\right) -\left( 5.4\right) $ show that 
\begin{equation}
u_{\varepsilon }\left( x,t\right) =e^{\varepsilon tW}u\left( x,t\right) 
\text{, for }t\in \left[ 0,1\right] .  \tag{5.5}
\end{equation}%
In particular, the equality $u_{\varepsilon }\left( x,1\right)
=e^{\varepsilon W}u\left( x,1\right) $, Lemma 3.1 with $a+ib=\varepsilon ,$ $%
\gamma =\frac{1}{\beta }$, $F\equiv 0$ and the fact that $u_{\varepsilon
}(0)=u\left( 0\right) $ imply that%
\[
\left\Vert e^{\frac{\left\vert x\right\vert ^{2}}{\beta ^{2}+4\varepsilon }%
}u_{\varepsilon }\left( .,1\right) \right\Vert _{X}\leq e^{\varepsilon
\left\Vert V_{1}\right\Vert _{B}}\left\Vert e^{\frac{\left\vert x\right\vert
^{2}}{\beta ^{2}}}u\left( .,1\right) \right\Vert _{X}\text{, }\left\Vert e^{%
\frac{\left\vert x\right\vert ^{2}}{^{\alpha ^{2}}}}u_{\varepsilon }\left(
.,0\right) \right\Vert _{X}=\left\Vert e^{\frac{\left\vert x\right\vert ^{2}%
}{^{\alpha ^{2}}}}u\left( .,0\right) \right\Vert _{X}\text{.} 
\]%
A second application of Lemma 3.1 with $a+ib$ $=$ $\varepsilon $, $F\equiv 0$%
, the value of $\gamma =\mu ^{-2}\left( t\right) $ and $\left( 5.2\right) $
show that%
\[
\left\Vert \varepsilon ^{\frac{\left\vert x\right\vert ^{2}}{\mu ^{2}\left(
t\right) +4\varepsilon t}}F_{\varepsilon }\left( .,t\right) \right\Vert
_{X}\leq e^{\varepsilon \left\Vert V_{1}\right\Vert _{B}}\left\Vert
\varepsilon ^{\frac{\left\vert x\right\vert ^{2}}{\mu ^{2}\left( t\right) }%
}V_{2}\left( .,t\right) \right\Vert _{B}\left\Vert u\left( .,t\right)
\right\Vert _{X}\text{, }t\in \left[ 0,1\right] . 
\]%
Setting, $\alpha _{\varepsilon }=\alpha +2\varepsilon $ and $\beta
_{\varepsilon }=\beta +2\varepsilon $, the last three inequalities give that%
\begin{equation}
\left\Vert e^{\frac{\left\vert x\right\vert ^{2}}{^{\beta _{\varepsilon
}^{2}}}}u_{\varepsilon }\left( .,1\right) \right\Vert _{X}\leq
e^{\varepsilon \left\Vert V_{1}\right\Vert _{B}}\left\Vert e^{\frac{%
\left\vert x\right\vert ^{2}}{^{\beta 2}}}u\left( .,1\right) \right\Vert
_{X},  \tag{5.6}
\end{equation}%
\[
\left\Vert e^{\frac{\left\vert x\right\vert ^{2}}{^{\alpha _{\varepsilon
}^{2}}}}u_{\varepsilon }\left( .,0\right) \right\Vert _{X}\leq
e^{\varepsilon \left\Vert V_{1}\right\Vert _{B}}\left\Vert e^{\frac{%
\left\vert x\right\vert ^{2}}{^{\alpha ^{2}}}}u\left( .,0\right) \right\Vert
_{X}, 
\]%
\begin{equation}
\left\Vert \varepsilon ^{\left\vert x\right\vert ^{2}\mu ^{-2}\left(
t\right) }F_{\varepsilon }\left( .,t\right) \right\Vert _{X}\leq
e^{\varepsilon \left\Vert V_{1}\right\Vert _{B}}\left\Vert \varepsilon
^{\left\vert x\right\vert ^{2}\mu ^{-2}\left( t\right) }V_{2}\left(
.,t\right) \right\Vert _{B}\left\Vert u\left( .,t\right) \right\Vert _{X}%
\text{, }t\in \left[ 0,1\right] .  \tag{5.7}
\end{equation}%
A third application of Lemma 3.1 with $a+ib=b,$ $F\equiv 0,$ $\gamma =0$,
and $(5.2),(5.5)$ implies that%
\begin{equation}
\left\Vert F_{\varepsilon }\left( .,t\right) \right\Vert _{X}\leq
e^{\varepsilon \left\Vert V_{1}\right\Vert _{B}}\left\Vert V_{2}\left(
.,t\right) \right\Vert _{B}\left\Vert u\left( .,t\right) \right\Vert _{X}%
\text{, }  \tag{5.8}
\end{equation}%
\[
\left\Vert u_{\varepsilon }\left( .,t\right) \right\Vert _{X}\leq
e^{\varepsilon \left\Vert V_{1}\right\Vert _{B}}\left\Vert u\left(
.,t\right) \right\Vert _{X}\text{, }t\in \left[ 0,1\right] . 
\]%
Set $\gamma _{\varepsilon }=\frac{1}{\alpha _{\varepsilon }\beta
_{\varepsilon }}$ and let%
\[
\tilde{u}_{\varepsilon }\left( x,t\right) =\left( \sqrt{\alpha _{\varepsilon
}\beta _{\varepsilon }}\rho _{\varepsilon }^{-1}\left( t\right) \right) ^{%
\frac{n}{2}}u\left( \sqrt{\alpha _{\varepsilon }\beta _{\varepsilon }}x\rho
_{\varepsilon }^{-1}\left( t\right) ,\beta _{\varepsilon }t\rho
_{\varepsilon }^{-1}\left( t\right) \right) e^{\varphi _{\varepsilon }} 
\]%
be the function associated to $u_{\varepsilon }$ in Lemma 4.1, where $%
a+ib=\varepsilon +i$ and $\alpha ,$ $\beta $ are replaced respectively by $%
\alpha _{\varepsilon }$, $\beta _{\varepsilon }$ when

\[
\rho _{\varepsilon }\left( t\right) =\alpha _{\varepsilon }\left( 1-t\right)
+\beta _{\varepsilon }t,\text{ }\varphi _{\varepsilon }=\varphi
_{\varepsilon }\left( x,t\right) =\frac{\left( \alpha _{\varepsilon }-\beta
_{\varepsilon }\right) \left\vert x\right\vert ^{2}}{4\left( a+ib\right)
\rho _{\varepsilon }\left( t\right) }. 
\]%
Because $\alpha <\beta $, $\tilde{u}_{\varepsilon }\in L^{\infty }\left(
0,1;X\right) \cap L^{2}\left( 0,1;Y^{1}\right) $ and satisfies the equation%
\[
\partial _{t}\tilde{u}_{\varepsilon }=\left( \varepsilon +i\right) \left(
\Delta \tilde{u}_{\varepsilon }+A\left( x\right) \tilde{u}_{\varepsilon }+%
\tilde{V}_{1}^{\varepsilon }\left( x,t\right) \tilde{u}_{\varepsilon }+%
\tilde{F}_{\varepsilon }\left( x,t\right) \right) \text{ in }R^{n}\times %
\left[ 0,1\right] , 
\]%
where 
\[
\tilde{V}_{1}^{\varepsilon }\left( x,t\right) =\alpha _{\varepsilon }\beta
_{\varepsilon }\rho _{\varepsilon }^{-2}\left( t\right) V_{1}\left( \sqrt{%
\alpha _{\varepsilon }\beta _{\varepsilon }}\rho _{\varepsilon }^{-1}\left(
t\right) x\right) ,\text{ }\sup\limits_{t\in \left[ 0,1\right] }\left\Vert 
\tilde{V}_{1}^{\varepsilon }\left( ,.t\right) \right\Vert _{B}\leq \frac{%
\beta }{\alpha }M_{1}, 
\]%
\begin{equation}
\tilde{F}_{\varepsilon }\left( x,t\right) =\left[ \sqrt{\alpha _{\varepsilon
}\beta _{\varepsilon }}\rho _{\varepsilon }^{-1}\left( t\right) \right] ^{%
\frac{n}{2}+2}F_{\varepsilon }\left( \sqrt{\alpha _{\varepsilon }\beta
_{\varepsilon }}\rho _{\varepsilon }^{-1}\left( t\right) x,\beta
_{\varepsilon }t\rho _{\varepsilon }^{-1}\left( t\right) \right) e^{\varphi
_{\varepsilon }},  \tag{5.9}
\end{equation}

\begin{equation}
\left\Vert e^{\gamma _{\varepsilon }\left\vert x\right\vert ^{2}}\tilde{F}%
_{\varepsilon }\left( .,t\right) \right\Vert _{X}\leq \frac{\beta }{\alpha }%
\left\Vert e^{\rho _{\varepsilon }^{-2}\left\vert x\right\vert
^{2}}F_{\varepsilon }\left( .,t\right) \right\Vert _{X}\text{, }\left\Vert 
\tilde{F}_{\varepsilon }\left( .,t\right) \right\Vert _{X}\leq \frac{\beta }{%
\alpha }\left\Vert F_{\varepsilon }\left( .,t\right) \right\Vert _{X} 
\tag{5.10}
\end{equation}%
and 
\begin{equation}
\left\Vert e^{\gamma _{\varepsilon }\left\vert x\right\vert ^{2}}\tilde{u}%
_{\varepsilon }\left( .,t\right) \right\Vert _{X}=\left\Vert e^{\left[ \rho
_{\varepsilon }^{-2}\left( s\right) +\varphi _{\varepsilon }\left(
s,t\right) \right] \left\vert y\right\vert ^{2}}u_{\varepsilon }\left(
.,s\right) \right\Vert _{X},\text{ }\left\Vert \tilde{u}_{\varepsilon
}\left( .,t\right) \right\Vert _{X}\leq \left\Vert u\left( .,s\right)
\right\Vert _{X},  \tag{5.11}
\end{equation}%
when $s=\beta _{\varepsilon }\rho _{\varepsilon }^{-1}\left( t\right) .$ The
above identity when $t$ is zero or one and $(5.6)$ shows that 
\begin{equation}
\left\Vert e^{\gamma _{\varepsilon }\left\vert x\right\vert ^{2}}\tilde{u}%
_{\varepsilon }\left( .,0\right) \right\Vert _{X}\leq \left\Vert e^{\frac{%
\left\vert x\right\vert ^{2}}{\beta ^{2}}}u\left( .,0\right) \right\Vert
_{X},\text{ }\left\Vert e^{\gamma _{\varepsilon }\left\vert x\right\vert
^{2}}\tilde{u}_{\varepsilon }\left( .,1\right) \right\Vert _{X}\leq 
\tag{5.12}
\end{equation}%
\[
e^{\varepsilon \left\Vert V\right\Vert _{B}}\left\Vert V_{2}\left(
.,t\right) \right\Vert _{B}\left\Vert e^{\frac{\left\vert x\right\vert ^{2}}{%
\beta ^{2}}}u\left( .,1\right) \right\Vert _{X}. 
\]%
On the other hand, 
\begin{equation}
N_{1}^{-1}\left\Vert u\left( .,0\right) \right\Vert _{X}\leq \left\Vert
u\left( .,t\right) \right\Vert _{X}\leq N_{1}\left\Vert u\left( .,0\right)
\right\Vert _{X},\text{ }t\in \left[ 0,1\right] ,  \tag{5.13}
\end{equation}%
where 
\[
N_{1}=e^{B\left( V_{2}\right) },\text{ }B\left( V_{2}\right)
=\sup\limits_{t\in \left[ 0,1\right] }\left\Vert \func{Re}V_{2}\left(
.,t\right) \right\Vert _{B}. 
\]%
The energy method imply that 
\begin{equation}
\partial _{t}\left\Vert \tilde{u}_{\varepsilon }\left( .,t\right)
\right\Vert _{X}^{2}\leq 2\varepsilon \left\Vert \tilde{V}_{1}^{\varepsilon
}\left( x,t\right) \right\Vert _{B}\left\Vert \tilde{u}_{\varepsilon }\left(
.,t\right) \right\Vert _{X}^{2}+2\left\Vert \tilde{F}_{\varepsilon }\left(
x,t\right) \right\Vert _{X}\left\Vert \tilde{u}_{\varepsilon }\left(
.,t\right) \right\Vert _{X}.  \tag{5.14}
\end{equation}

Let $0=t_{0}<t_{1}<....t_{m}=1$ be a uniformly distributed partition of $%
[0,1]$, where $m$ will be chosen later. The inequalities $\left( 5.8\right)
-\left( 5.10\right) ,$ $\left( 5.13\right) $ and $(5.14$) imply that there
is $N_{2}$, which depends on $\frac{\beta }{\alpha }$, $\left\Vert
V_{1}\right\Vert _{B}$ and $B\left( V_{2}\right) $ such that%
\begin{equation}
\left\Vert \tilde{u}_{\varepsilon }\left( .,t_{i}\right) \right\Vert
_{X}\leq e^{\frac{\varepsilon \beta }{\alpha }\left\Vert V_{1}\right\Vert
_{B}}\left\Vert \tilde{u}_{\varepsilon }\left( .,t\right) \right\Vert
_{X}+N_{2}\sqrt{t_{i}-t_{i-1}}\left\Vert u\left( .,0\right) \right\Vert _{X}
\tag{5.15}
\end{equation}%
for $t\in \left[ t_{i-1},t_{i}\right] $ and $i=1,2,...m.$ Choose now $m$ so
that%
\begin{equation}
N_{2}\max\limits_{i}\sqrt{t_{i}-t_{i-1}}\leq \frac{1}{4N_{1}}.  \tag{5.16}
\end{equation}%
Because, $\lim\limits_{\varepsilon \rightarrow 0}\left\Vert \tilde{u}%
_{\varepsilon }\left( .,t\right) \right\Vert _{X}=\left\Vert u\left(
.,s\right) \right\Vert _{X}$ when $s=\beta t\rho \left( t\right) $ and $%
(5.13)$, there is $\varepsilon _{0}$ such that 
\begin{equation}
\left\Vert \tilde{u}_{\varepsilon }\left( .,t_{i}\right) \right\Vert
_{X}\geq \frac{1}{4N_{1}}\left\Vert u\left( .,0\right) \right\Vert _{X},%
\text{ when }0<\varepsilon \leq \varepsilon _{0},\text{ }i=1,2,...m 
\tag{5.17}
\end{equation}%
and now, $(5.15)-(5.17)$ show that 
\begin{equation}
\left\Vert \tilde{u}_{\varepsilon }\left( .,t\right) \right\Vert _{X}\geq 
\frac{1}{4N_{1}}\left\Vert u\left( .,0\right) \right\Vert _{X},\text{ when }%
0<\varepsilon \leq \varepsilon _{0},\text{ }t\in \left[ 0,1\right] . 
\tag{5.18}
\end{equation}%
It is now simple to verify that $(5.18)$, the first inequality in $(5.7),$ $%
(5.10)$ and $(5.13)$ imply that%
\begin{equation}
\sup\limits_{t\in \left[ 0,1\right] }\frac{\left\Vert e^{\gamma
_{\varepsilon }\left\vert x\right\vert ^{2}}\tilde{F}_{\varepsilon }\left(
.,t\right) \right\Vert _{X}}{\left\Vert \tilde{u}_{\varepsilon }\left(
.,t\right) \right\Vert _{X}}\leq \frac{4\beta }{\alpha }M_{2}\left(
\varepsilon \right) ,\text{ when }0<\varepsilon \leq \varepsilon _{0}, 
\tag{5.19}
\end{equation}%
where 
\[
M_{2}\left( \varepsilon \right) =e^{\sup\limits_{t\in \left[ 0,1\right]
}\left\Vert \func{Re}V_{2}\left( .,t\right) \right\Vert _{B}+\varepsilon
\left\Vert V_{1}\right\Vert _{B}}\sup\limits_{t\in \left[ 0,1\right]
}\left\Vert e^{\left\vert x\right\vert ^{2}\mu ^{-2}\left( t\right)
}V_{2}\left( .,t\right) \right\Vert _{B}. 
\]%
By using Lemma $3$.$3,$ $(5.12),$ $(5.9)$ and $(5.19)$ to show that $%
\left\Vert e^{\gamma _{\varepsilon }\left\vert x\right\vert ^{2}}\tilde{u}%
_{\varepsilon }\left( .,t\right) \right\Vert _{X}$ is logarithmically convex
in $[0,1]$ and that%
\begin{equation}
\left\Vert e^{\gamma \left\vert x\right\vert ^{2}}\tilde{u}_{\varepsilon
}\left( .,t\right) \right\Vert _{X}\leq e^{NM\left( a,b\right) }\left\Vert
e^{\gamma \left\vert x\right\vert ^{2}}\tilde{u}_{\varepsilon }\left(
0\right) \right\Vert _{X}^{1-t}\left\Vert e^{\gamma \left\vert x\right\vert
^{2}}\tilde{u}_{\varepsilon }\left( 1\right) \right\Vert _{X}^{t}, 
\tag{5.20}
\end{equation}%
when $0<\varepsilon \leq \varepsilon _{0},$ $t\in \left[ 0,1\right] $ and $%
N=N\left( \alpha ,\beta \right) .$ Then, Lemma 3.4 gives that%
\[
\left\Vert \eta \nabla \tilde{u}_{\varepsilon }\right\Vert _{Z}+\left\Vert
\eta \left\vert x\right\vert \tilde{u}_{\varepsilon }\right\Vert _{Z}\leq 
\]

\[
N\left( 1+M_{1}\right) \left[ \sup\limits_{t\in \left[ 0,1\right]
}\left\Vert e^{\gamma \left\vert x\right\vert ^{2}}\tilde{u}_{\varepsilon
}\left( .,t\right) \right\Vert _{X}+\sup\limits_{t\in \left[ 0,1\right]
}\left\Vert e^{\gamma \left\vert x\right\vert ^{2}}\tilde{F}_{\varepsilon
}\left( .,t\right) \right\Vert _{Z}\right] \leq 
\]%
\[
Ne^{N\left( M_{0}+M_{1}+M_{2}\left( \varepsilon \right)
+M_{1}^{2}+M_{2}^{2}\left( \varepsilon \right) \right) }\left[ \left\Vert e^{%
\frac{\left\vert x\right\vert ^{2}}{\beta ^{2}}}u\left( .,0\right)
\right\Vert _{X}+\left\Vert e^{\frac{\left\vert x\right\vert ^{2}}{\alpha
^{2}}}u\left( .,1\right) \right\Vert _{X}\right] , 
\]%
when $0<\varepsilon \leq \varepsilon _{0}$, the logarithmic convexity and
regularity of $u$ follow from the limit of the identity in $(5.11)$, the
final limit relation between the variables $s$ and $t$, $s=$ $\beta t\rho
\left( t\right) $ and letting $\varepsilon $ tend to zero in $(5.20)$ and
the above inequality.

By reasoning as in $\left[ \text{4, Lemma 6}\right] $ we obtain:

\textbf{Lemma 5.1.} Let $A=A\left( x\right) $ be a symmetric operator in
Hilbert space $H$ with independent on $x\in R^{n}$ domain $D\left( A\right) $
that is dense on $H$ and 
\[
\left\Vert V\right\Vert _{B}\leq \varepsilon _{0}\text{ for a }\varepsilon
_{0}>0. 
\]

Let $u\in C\left( \left[ 0,1\right] ;X\left( A\right) \right) $ be a
solution of the equation 
\[
\partial _{t}u=i\left[ \Delta u+Au+V\left( x,t\right) u+F\left( x,t\right) %
\right] ,\text{ }x\in R^{n},\text{ }t\in \left[ 0,1\right] . 
\]

\bigskip Then, 
\[
\sup\limits_{t\in \left[ 0,1\right] }\left\Vert e^{\lambda .x}u\left(
.,t\right) \right\Vert _{X}\leq N\left[ \left\Vert e^{\lambda .x}u\left(
.,0\right) \right\Vert _{X}+\left\Vert e^{\lambda .x}u\left( .,1\right)
\right\Vert _{X}+\left\Vert e^{\lambda .x}F\left( .,t\right) \right\Vert
_{L^{1}\left( 0,1;X\right) }\right] , 
\]%
where $\lambda \in R^{n}$ and $N>0$ is a constant.

\textbf{Theorem 5.1.} Assume that $A=A\left( x\right) $ is a symmetric
operator in Hilbert space $H$ with independent on $x\in R^{n}$ domain $%
D\left( A\right) $ that is dense on $H$ and 
\[
V\in B\text{ and }\lim\limits_{r\rightarrow \infty }\left\Vert V\right\Vert
_{O\left( r\right) }=0. 
\]

Suppose that $\alpha ,$ $\beta $ are positive numbers and 
\[
\left\Vert e^{\frac{\left\vert x\right\vert ^{2}}{\beta ^{2}}}u\left(
.,0\right) \right\Vert _{X}<\infty \text{, }\left\Vert e^{\frac{\left\vert
x\right\vert ^{2}}{\alpha ^{2}}}u\left( .,1\right) \right\Vert _{X}<\infty . 
\]

Let $u\in C\left( \left[ 0,1\right] ;X\left( A\right) \right) $ be a
solution of the equation 
\[
\partial _{t}u=i\left[ \Delta u+A\left( x\right) u+V\left( x,t\right) u%
\right] ,\text{ }x\in R^{n},\text{ }t\in \left[ 0,1\right] . 
\]

\bigskip Then, there is a $N=N(\alpha ,\beta )$ such that 
\[
\sup\limits_{t\in \left[ 0,1\right] }\left\Vert e^{\left\vert x\right\vert
^{2}\mu ^{-2}\left( t\right) }u\left( .,t\right) \right\Vert _{X}+\left\Vert 
\sqrt{t\left( 1-t\right) }e^{\left\vert x\right\vert ^{2}\mu ^{-2}\left(
t\right) }\nabla u\right\Vert _{L^{2}\left( R^{n}\times \left[ 0,1\right]
;H\right) }\leq 
\]%
\[
Ne^{B\left( V\right) }\left[ \left\Vert e^{\frac{\left\vert x\right\vert ^{2}%
}{\beta ^{2}}}u\left( .,0\right) \right\Vert _{X}+\left\Vert e^{\frac{%
\left\vert x\right\vert ^{2}}{\alpha ^{2}}}u\left( .,1\right) \right\Vert
_{X}+\sup\limits_{t\in \left[ 0,1\right] }\left\Vert u\left( .,t\right)
\right\Vert _{X}\right] ,\text{ } 
\]%
where, 
\[
B\left( V\right) =\sup\limits_{t\in \left[ 0,1\right] }\left\Vert
V\right\Vert _{B}. 
\]

\textbf{Proof. }Assume that $u(y,s)$ verifies the equation%
\[
\partial _{s}u=i\left[ \Delta u+A\left( y\right) u+V\left( y,s\right)
u+F\left( y,s\right) \right] ,\text{ }y\in R^{n},\text{ }s\in \left[ 0,1%
\right] . 
\]%
Set $\gamma =\left( \alpha \beta \right) ^{-1}$ and let 
\begin{equation}
\tilde{u}\left( x,t\right) =\left( \sqrt{\alpha \beta }\rho ^{-1}\left(
t\right) \right) ^{\frac{n}{2}}u\left( \sqrt{\alpha \beta }x\rho ^{-1}\left(
t\right) ,\beta t\rho ^{-1}\left( t\right) \right) e^{\varphi }.  \tag{5.21}
\end{equation}%
The function $\left( 5.21\right) $ is a solution of 
\[
\partial _{t}u=i\left[ \Delta u+A\left( x\right) u+V\left( x,t\right) u%
\right] ,\text{ }x\in R^{n},\text{ }t\in \left[ 0,1\right] 
\]%
with 
\[
\tilde{V}\left( x,t\right) =\alpha \beta \delta ^{-2}\left( t\right) V\left( 
\sqrt{\alpha \beta }x\rho ^{-1}\left( t\right) ,\beta t\rho ^{-1}\left(
t\right) \right) , 
\]%
\[
\sup\limits_{t\in \left[ 0,1\right] }\left\Vert \tilde{V}\left( .,t\right)
\right\Vert _{B}\leq \max \left( \frac{\alpha }{\beta },\frac{\beta }{\alpha 
}\right) \sup\limits_{t\in \left[ 0,1\right] }\left\Vert V\left( .,t\right)
\right\Vert _{B},\text{ }\lim\limits_{r\rightarrow \infty }\left\Vert \tilde{%
V}\left( .,t\right) \right\Vert _{O\left( r\right) }=0 
\]%
and%
\begin{equation}
\left\Vert e^{\gamma \left\vert x\right\vert ^{2}}\tilde{u}\left( .,t\right)
\right\Vert _{X}=\left\Vert e^{\mu ^{2}\left( t\right) \left\vert
x\right\vert ^{2}}u\left( .,s\right) \right\Vert _{X},  \tag{5.22 }
\end{equation}%
\[
\left\Vert \tilde{u}\left( .,t\right) \right\Vert _{X}=\left\Vert u\left(
.,s\right) \right\Vert _{X}\text{ when }s=\beta t\mu \left( t\right) . 
\]%
Choose $r>0$ such that $\left\Vert \tilde{V}\left( .,t\right) \right\Vert
_{O\left( r\right) }\leq \varepsilon _{0}$ we get%
\[
\partial _{t}\tilde{u}=i\left[ \Delta \tilde{u}+A\tilde{u}+\tilde{V}%
_{r}\left( x,t\right) u+\tilde{F}_{r}\left( x,t\right) \right] ,\text{ }x\in
R^{n},\text{ }t\in \left[ 0,1\right] , 
\]%
with 
\[
\tilde{V}_{r}\left( x,t\right) =\chi _{R^{n}/O_{r}}\tilde{V}\left(
x,t\right) \text{, }\tilde{F}_{r}\left( x,t\right) =\chi _{O_{r}}\tilde{V}%
\left( x,t\right) \tilde{u}. 
\]%
Then using the Lemma 5.1 we obtain 
\[
\sup\limits_{t\in \left[ 0,1\right] }\left\Vert e^{\lambda .x}\tilde{u}%
\left( .,t\right) \right\Vert _{X}\leq 
\]%
\[
N\left[ \left\Vert e^{\lambda .x}\tilde{u}\left( .,0\right) \right\Vert
_{X}+\left\Vert e^{\lambda .x}\tilde{u}\left( .,1\right) \right\Vert
_{X}+e^{\left\vert \lambda \right\vert r}\left\Vert \tilde{V}\left(
.,t\right) \right\Vert _{B}\sup\limits_{t\in \left[ 0,1\right] }\left\Vert
u\left( .,t\right) \right\Vert _{X}\right] . 
\]%
Replace $\lambda $ by $\lambda \sqrt{\gamma }$ in the above inequality,
square both sides, multiply all by $e^{-\frac{\left\vert \lambda
^{2}\right\vert }{2}}$ and integrate both sides with respect to $\lambda $
in $R^{n}$. This and the identity%
\[
\dint\limits_{R^{n}}e^{2\sqrt{\gamma }\lambda .x-\frac{\left\vert \lambda
\right\vert ^{2}}{2}}d\lambda =\left( 2\pi \right) ^{\frac{n}{2}}e^{2\gamma
\left\vert x\right\vert ^{2}} 
\]%
imply the inequality 
\begin{equation}
\sup\limits_{t\in \left[ 0,1\right] }\left\Vert \tilde{u}\left( .,t\right)
\right\Vert _{X}\leq  \tag{5.23}
\end{equation}%
\[
N\left[ \left\Vert e^{2\gamma \left\vert x\right\vert ^{2}}\tilde{u}\left(
.,0\right) \right\Vert _{X}+\left\Vert e^{2\gamma \left\vert x\right\vert
^{2}}\tilde{u}\left( .,1\right) \right\Vert _{X}+\left\Vert e^{2\gamma r^{2}}%
\tilde{V}\left( .,t\right) \right\Vert _{B}\sup\limits_{t\in \left[ 0,1%
\right] }\left\Vert \tilde{u}\left( .,t\right) \right\Vert _{X}\right] . 
\]%
This inequality and $(5.22)$ imply that%
\[
\sup\limits_{t\in \left[ 0,1\right] }\left\Vert \tilde{u}\left( .,t\right)
\right\Vert _{X}\leq 
\]%
\[
N\left[ \left\Vert e^{\frac{\left\vert x\right\vert ^{2}}{\beta ^{2}}}\tilde{%
u}\left( .,0\right) \right\Vert _{X}+\left\Vert e^{\frac{\left\vert
x\right\vert ^{2}}{\beta ^{2}}}\tilde{u}\left( .,1\right) \right\Vert
_{X}+\sup\limits_{t\in \left[ 0,1\right] }\left\Vert V\left( .,t\right)
\right\Vert _{B}\sup\limits_{t\in \left[ 0,1\right] }\left\Vert u\left(
.,t\right) \right\Vert _{X}\right] 
\]%
for some new constant $N$.

To prove the regularity of $u$ we proceed as in $(5.2)-$ $(5.4)$. The
Duhamel formula shows that%
\begin{equation}
u_{\varepsilon }\left( x,t\right)
=e^{itW}u_{0}+i\dint\limits_{0}^{t}e^{i\left( t-s\right) W}V_{2}\left(
x,s\right) u\left( x,s\right) ds,\text{ }x\in R^{n},\text{ }t\in \left[ 0,1%
\right] .  \tag{5.24}
\end{equation}%
For $0\leq \varepsilon \leq 1$, set 
\begin{equation}
\tilde{F}_{\varepsilon }\left( x,t\right) =\frac{i}{\varepsilon +i}%
e^{\varepsilon t\left( \Delta +A\right) }\tilde{V}\left( x,t\right) \tilde{u}%
\left( x,t\right) ,  \tag{5.25}
\end{equation}%
and 
\begin{equation}
\tilde{u}_{\varepsilon }\left( x,t\right) =e^{\left( \varepsilon +i\right)
t\left( \Delta +A\right) }u_{0}+\left( \varepsilon +i\right)
\dint\limits_{0}^{t}e^{\left( \varepsilon +i\right) \left( t-s\right) \left(
\Delta +A\right) }\tilde{F}\left( x,s\right) u\left( x,s\right) ds, 
\tag{5.26}
\end{equation}%
\[
\text{ }x\in R^{n},\text{ }t\in \left[ 0,1\right] . 
\]%
The identities 
\[
e^{\left( z_{1}+z_{2}\right) \left( \Delta +A\right) }=e^{z_{1}\left( \Delta
+A\right) }.e^{z_{2}\left( \Delta +A\right) }\text{ for }\func{Re}z_{1},%
\text{ }\func{Re}z_{2}\geq 0 
\]%
and $\left( 5.24\right) -\left( 5.26\right) $ show that%
\begin{equation}
\tilde{u}_{\varepsilon }\left( x,t\right) =e^{\varepsilon t\left( \Delta
+A\right) }\tilde{u}\left( x,t\right) \text{ for }t\in \left[ 0,1\right] . 
\tag{5.27}
\end{equation}%
From Lemma 3.1 with $a+ib=\varepsilon $, $\left( 5.27\right) $ and $\left(
5.25\right) $ we get that 
\begin{equation}
\sup\limits_{t\in \left[ 0,1\right] }\left\Vert e^{\gamma _{\varepsilon
}\left\vert x\right\vert ^{2}}\tilde{u}_{\varepsilon }\left( .,t\right)
\right\Vert _{X}\leq \sup\limits_{t\in \left[ 0,1\right] }\left\Vert
e^{\gamma \left\vert x\right\vert ^{2}}\tilde{u}\left( .,t\right)
\right\Vert _{X},  \tag{4.28}
\end{equation}%
\[
\sup\limits_{t\in \left[ 0,1\right] }\left\Vert e^{\gamma _{\varepsilon
}\left\vert x\right\vert ^{2}}\tilde{F}_{\varepsilon }\left( .,t\right)
\right\Vert _{X}\leq e^{\tilde{V}_{0}}\sup\limits_{t\in \left[ 0,1\right]
}\left\Vert e^{\gamma \left\vert x\right\vert ^{2}}\tilde{F}\left(
.,t\right) \right\Vert _{X}, 
\]%
where 
\[
\gamma _{\varepsilon }=\frac{\gamma }{1+4\gamma \varepsilon },\text{ }\tilde{%
V}_{0}=\sup\limits_{t\in \left[ 0,1\right] }\left\Vert \tilde{V}\right\Vert
_{B}. 
\]%
Then, Lemma 3.4, $(5.28)$ and $(5.23)$ show that 
\[
\left\Vert e^{\gamma _{\varepsilon }\left\vert x\right\vert ^{2}}u\left(
.,t\right) \right\Vert _{L^{2}\left( R^{n}\times \left[ 0,1\right] ;H\right)
}+\left\Vert \sqrt{t\left( 1-t\right) }e^{\gamma _{\varepsilon }\left\vert
x\right\vert ^{2}}\nabla u\right\Vert _{L^{2}\left( R^{n}\times \left[ 0,1%
\right] ;H\right) }\leq 
\]%
\[
Ne^{NV_{0}}\left[ \left\Vert e^{\frac{\left\vert x\right\vert ^{2}}{\beta
^{2}}}u\left( .,0\right) \right\Vert _{X}+\left\Vert e^{\frac{\left\vert
x\right\vert ^{2}}{\alpha ^{2}}}u\left( .,1\right) \right\Vert
_{X}+\sup\limits_{t\in \left[ 0,1\right] }\left\Vert u\left( .,t\right)
\right\Vert _{X}\right] , 
\]%
where 
\[
V_{0}=\sup\limits_{t\in \left[ 0,1\right] }\left\Vert V\left( x,t\right)
\right\Vert _{B}. 
\]%
The Theorem 5.1 follows from this inequality, from $(5.21)-(5.23)$ and
letting $\varepsilon $ tend to zero.

\begin{center}
\textbf{6. A Hardy type abstract uncertainty principle. Proof of Theorem 1.}
\end{center}

\bigskip The assertion about the Carleman inequality in Lemma 6.1 below is
the following monotonicity or frequency function argument related to Lemma
3. 2. When $u\in C([0,1];X)$ is a solution to the free abstract Schr\"{o}%
dinger equation%
\[
\partial _{t}u-i\left( \Delta u+A\left( x\right) u\right) =0,\text{ }x\in
R^{n},\text{ }t\in \left[ 0,1\right] , 
\]%
satisfies

\[
\left\Vert e^{\gamma \left\vert x\right\vert ^{2}}u\left( .,0\right)
\right\Vert _{X}<\infty \text{, }\left\Vert e^{\gamma \left\vert
x\right\vert ^{2}}u\left( .,1\right) \right\Vert _{X}<\infty 
\]%
and 
\[
\ f=e^{\varkappa }u,\text{ }Q\left( t\right) =\left( f\left( .,t\right)
,f\left( .,t\right) \right) _{X}, 
\]%
where

\[
\varkappa \left( x,t\right) =%
\mu
|x+rt(1-t)|^{2}-\frac{r^{2}t(1-t)}{8\mu },\text{ }\sigma \left( \varepsilon
,t\right) =\frac{\left( 1+\varepsilon \right) t(1-t)}{16\mu }. 
\]%
Then, $\log Q\left( t\right) $ is logaritmicaly convex in $[0,1],$ when $%
0<\mu <\gamma .$

The formal application of the above argument to a $C([0,1];X)$ solution of
the equation%
\begin{equation}
\partial _{t}u-i\left[ \Delta u+Au+V\left( x,t\right) u\right] =0,\text{ }%
x\in R^{n},\text{ }t\in \left[ 0,1\right] ,  \tag{6.1}
\end{equation}%
implies a similar result, when $V$ is a bounded potential, though the
justification of the correctness of the assertions involved in the
corresponding formal application of Lemma 3.2 were formal. In fact, we can
only justify these assertions, when the potential $V$ verifies the first
condition in Theorem 1 or when we can obtain the additional regularity of
the gradient of $u$ in the strip, as in Theorem 5.1. Here, we choose to
prove Theorem 1 using the Carleman inequality in Lemma 6.1 in place of the
above convexity argument. The reason for our choice is that it is simpler to
justify the correctness of the application of the Carleman inequality to a $%
C([0,1];X)$ solution to $(6.1)$ than the corresponding monotonicity or
logarithmic convexity of the solution.

\textbf{Lemma 6.1. }Let the assumptions (1)-(2) of Conditon 3.1 hold.
Moreover,

\[
V\in B\text{ and }\lim\limits_{r\rightarrow \infty }\left\Vert V\right\Vert
_{O\left( r\right) }=0. 
\]%
Then the estimate

\[
r\sqrt{\frac{\varepsilon }{8\mu }}\left\Vert e^{\varkappa -\sigma }\upsilon
\right\Vert _{L^{2}\left( R^{n+1};H\right) }\leq \left\Vert e^{\varkappa
-\sigma }\left[ \partial _{t}u-i\left( \Delta u+Au\right) \right] \upsilon
\right\Vert _{L^{2}\left( R^{n+1};H\right) } 
\]%
holds, when $\varepsilon >0$, $\mu >0$, $r>0$ and $\upsilon \in
C_{0}^{\infty }\left( R^{n+1};H\left( A\right) \right) $.

\textbf{Proof.} Let\textbf{\ }$f=e^{\varkappa -\sigma }\upsilon $. Then, 
\[
e^{\varkappa -\sigma }\left[ \partial _{t}u-i\left( \Delta u+Au\right) %
\right] \upsilon =\partial _{t}f+Sf-Kf. 
\]%
From $(3.8)-(3.10)$ with $\gamma =1$, $a+ib=i$ and $\varphi \left(
x,t\right) =\varkappa \left( x,t\right) -\sigma \left( \varepsilon ,t\right) 
$ we have 
\[
S=-4%
\mu
i(x+rt(1-t)e_{1})\text{\textperiodcentered }\nabla -2%
\mu
ni+2%
\mu
r(1-2t)(x_{1}+rt(1-t))-\sigma ,\text{ } 
\]%
\[
K=i\left( \triangle +A\right) +4%
\mu
^{2}i|x+rt(1-t)e_{1}|^{2}\text{, }S_{t}+[S,K]= 
\]%
\[
-8%
\mu
\triangle +32%
\mu
^{3}|x+rt(1-t)e_{1}|^{2}-4%
\mu
r(x_{1}+rt(1-t))+ 
\]

\[
+2%
\mu
r^{2}(1-2t)^{2}+\frac{\left( 1+\varepsilon \right) r^{2}}{8\mu }+-4i%
\mu
r(1-2t)\partial _{x_{1}} 
\]%
and

\begin{equation}
(S_{t}f+[S,K]f,f)_{X}=32%
\mu
^{3}\dint\limits_{R^{n}}\left\vert x+rt(1-t)e_{1}-\frac{r}{16\mu ^{2}}%
e_{1}\right\vert ^{2}\left\Vert f\right\Vert ^{2}dx+  \tag{6.2}
\end{equation}%
\[
\frac{\varepsilon r^{2}}{8\mu }\dint\limits_{R^{n}}\left\Vert f\right\Vert
^{2}dx+8\mu \dint\limits_{R^{n}}\left\Vert \nabla _{x^{\prime }}f\right\Vert
_{H}^{2}dx+8\mu \dint\limits_{R^{n}}\left\Vert i\partial _{x_{1}}f-r\left( 
\frac{1}{2}-t\right) f\right\Vert ^{2}dx\geq 
\]%
\[
\frac{\varepsilon r^{2}}{8\mu }\dint\limits_{R^{n}}\left\Vert f\right\Vert
^{2}dx. 
\]%
Following the standard method to handle $L_{2}$-Carleman inequalities, the
symmetric and skew-symmetric parts of $\partial _{t}-S-K$, as a space-time
operator, are respectively $-S$ and $\partial _{t}-K$, and $[-S,\partial
_{t}-K]=S_{t}+[S,K]$. Thus, 
\[
\left\Vert \partial _{t}f-Sf-Kf\right\Vert _{L^{2}\left( R^{n+1};H\right)
}^{2}=\left\Vert \partial _{t}f-Kf\right\Vert _{L^{2}\left( R^{n+1};H\right)
}^{2}+\left\Vert Sf\right\Vert _{L^{2}\left( R^{n+1};H\right) }^{2}- 
\]%
\begin{equation}
2\func{Re}\dint\limits_{R^{n}}\dint\limits_{-\infty }^{\infty }\left(
Sf,\partial _{t}f-Kf\right) dxdt\geq
\dint\limits_{R^{n}}\dint\limits_{-\infty }^{\infty }\left( [-S,\partial
_{t}-K]f,f\right) dxdt=  \tag{6.3}
\end{equation}%
\[
\dint\limits_{-\infty }^{\infty }\left( S_{t}f+\left[ S,K\right] f,f\right)
_{H}dt, 
\]%
and the Lemma 6.1 follows from $(5.2)$ and $(5.3)$.

\textbf{Proof of Theorem 1. }Let $u$ be as in Theorem 1, $\tilde{u}$ and $%
\tilde{V}$ be corresponding functions defined in Lemma 4.1, when $a+ib=i$.
Then, $\tilde{u}\in $ $C([0,1];X\left( A\right) )$ is a solution of the
equation%
\[
\partial _{t}u-i\left[ \Delta u+Au+\tilde{V}u\right] =0,\text{ }x\in R^{n},%
\text{ }t\in \left[ 0,1\right] 
\]%
and 
\[
\left\Vert e^{\gamma \left\vert x\right\vert ^{2}}\tilde{u}\left( .,0\right)
\right\Vert _{X}<\infty \text{, }\left\Vert e^{\gamma \left\vert
x\right\vert ^{2}}\tilde{u}\left( .,1\right) \right\Vert _{X}<\infty \text{
for }\gamma =\frac{1}{\alpha \beta },\text{ }\gamma >\frac{1}{2}. 
\]%
The proof of Theorem 3 show that in either case%
\begin{equation}
N_{\gamma }=\sup\limits_{t\in \left[ 0,1\right] }\left[ \left\Vert e^{\gamma
_{\varepsilon }\left\vert x\right\vert ^{2}}\tilde{u}\left( .,t\right)
\right\Vert _{L^{2}\left( R^{n}\times \left[ 0,1\right] ;H\right)
}+\left\Vert \sqrt{t\left( 1-t\right) }e^{\gamma _{\varepsilon }\left\vert
x\right\vert ^{2}}\nabla \tilde{u}\right\Vert _{L^{2}\left( R^{n}\times %
\left[ 0,1\right] ;H\right) }\right] <\infty .  \tag{6.4}
\end{equation}%
For given $r>0$, choose $%
\mu
$ and $\varepsilon $ such that%
\begin{equation}
\frac{\left( 1+\varepsilon \right) ^{\frac{3}{2}}}{2\left( 1-\varepsilon
\right) ^{3}}<\mu \leq \frac{\gamma }{1+\varepsilon }  \tag{6.5}
\end{equation}%
and let $\eta _{M}$ and $\theta _{r}$ be smooth functions verifying, $\theta
_{M}\left( x\right) $ $=1$, when $|x|\leq M$, $\theta _{M}\left( x\right) $ $%
=0$, when $|x|>2M$, $M\geq 2r$, $\eta _{r}\in C_{0}^{\infty }(0,1),$ $0\leq
\eta _{r}\left( t\right) \leq 1$, $\eta _{r}\left( t\right) =1$ for $t\in %
\left[ \frac{1}{r},1-\frac{1}{r}\right] $ and $\eta _{r}\left( t\right) =0$
for $t\in \left[ 0,\frac{1}{2r}\right] \cup \left[ 1-\frac{1}{2r},1\right] .$
Then, $\upsilon \left( x,t\right) =\eta _{r}\left( t\right) \theta
_{M}\left( x\right) \tilde{u}\left( x,t\right) $ is compactly supported in $%
R^{n}\times (0,1)$ and 
\begin{equation}
\partial _{t}\upsilon -i\left[ \Delta \upsilon +A\upsilon +\tilde{V}\upsilon %
\right] =\eta _{r}^{\prime }\left( t\right) \theta _{M}\left( x\right) 
\tilde{u}\left( x,t\right) -\left( 2\nabla \theta _{M}.\nabla \tilde{u}+%
\tilde{u}\Delta \theta _{M}\right) \eta _{r}.\text{ }  \tag{6.6}
\end{equation}%
The terms on the right hand side of $(6.6)$ are supported, where 
\[
\mu
|x+rt(1-t)|^{2}\leq \gamma \left\vert x\right\vert ^{2}+\frac{\gamma }{%
\varepsilon }, 
\]%
\[
\mu
|x+rt(1-t)e_{1}|^{2}\leq \gamma \left\vert x\right\vert ^{2}+\frac{\gamma }{%
\varepsilon }r^{2}. 
\]%
Apply now Lemma 6.1 to $\upsilon $ with the values of $%
\mu
$ and $\varepsilon $ chosen in $(6.5)$. This, the bounds for $%
\mu
|x+rt(1-t)e_{1}|^{2}$ in each of the parts of the support of%
\[
\partial _{t}\upsilon -i\left[ \Delta \upsilon +A\upsilon +\tilde{V}\upsilon %
\right] 
\]%
and the natural bounds for $\nabla \theta _{M}$, $\triangle \theta _{M}$ and 
$\eta _{r}^{\prime }$ show that there is a constant $N_{\varepsilon }$ such
that 
\[
r\left\Vert e^{\varkappa -\sigma }\upsilon \right\Vert _{L^{\infty }\left(
R^{n}\times \left[ 0,1\right] ;H\right) }\leq 
\]%
\begin{equation}
N_{\varepsilon }\left\Vert \tilde{V}\right\Vert _{B}\left\Vert e^{\varkappa
-\sigma }\upsilon \right\Vert _{L^{2}\left( R^{n}\times \left[ 0,1\right]
;H\right) }+N_{\varepsilon }re^{\frac{\gamma }{\varepsilon }%
}\sup\limits_{t\in \left[ 0,1\right] }\left\Vert e^{\gamma \left\vert
x\right\vert ^{2}}\tilde{u}\left( .,t\right) \right\Vert _{X}+  \tag{6.7}
\end{equation}%
\[
N_{\varepsilon }M^{-1}e^{\frac{\gamma }{\varepsilon }r^{2}}\left\Vert
e^{\gamma \left\vert x\right\vert ^{2}}\left( \left\Vert \tilde{u}%
\right\Vert +\left\Vert \nabla \tilde{u}\right\Vert _{H}\right) \right\Vert
_{L^{2}\left( R^{n}\times \sigma _{r}\right) }, 
\]%
where 
\[
\sigma _{r}=\left[ \frac{1}{2r},1-\frac{1}{2r}\right] . 
\]%
The first term on the right hand side of $(6.7)$ can be hidden in the left
hand side, when $r\geq 2N_{\varepsilon }\left\Vert \tilde{V}\right\Vert _{B}$%
, while the last tends to zero, when $M$ tends to infinity by $(6.4)$. This
and the fact that $\upsilon =\tilde{u}$ in $O_{r_{\varepsilon }}\times \left[
\frac{1-\varepsilon }{2},\frac{1+\varepsilon }{2}\right] ,$ where 
\[
\varkappa -\sigma \geq \frac{r^{2}}{16\mu }\left( 4\mu ^{2}\left(
1-\varepsilon \right) ^{6}-\left( 1+\varepsilon \right) ^{3}\right) ,\text{ }%
r_{\varepsilon }=\frac{\varepsilon \left( 1-\varepsilon ^{2}\right) ^{2}r}{4}%
; 
\]%
and $\left( 6.5\right) $ show that%
\begin{equation}
e^{C\left( \gamma ,\varepsilon \right) }\left\Vert \tilde{u}\right\Vert
_{L^{2}\left( O\left( r,\varepsilon \right) ;H\right) }\leq N_{\gamma
,\varepsilon },\text{ }O\left( r,\varepsilon \right) =O_{\frac{r}{8}}\times %
\left[ \frac{1-\varepsilon }{2},\frac{1+\varepsilon }{2}\right]  \tag{6.8}
\end{equation}%
when $r\geq 2N_{\varepsilon }\left\Vert \tilde{V}\right\Vert _{B}.$ At the
same time 
\begin{equation}
\left( B\left( \tilde{V}\right) \right) ^{-1}\left\Vert \tilde{u}\left(
.,0\right) \right\Vert _{X}\leq \left\Vert \tilde{u}\left( .,t\right)
\right\Vert _{X}\leq B\left( \tilde{V}\right) \left\Vert \tilde{u}\left(
.,1\right) \right\Vert _{X}  \tag{6.9}
\end{equation}%
for $0\leq t\leq 1$ and $B\left( \tilde{V}\right) =\sup\limits_{t\in \left[
0,1\right] }\left\Vert \tilde{V}\right\Vert _{B}.$ Moreover, from $\left(
6.4\right) $ we get 
\begin{equation}
\left\Vert \tilde{u}\left( .,t\right) \right\Vert _{X}\leq \left\Vert \tilde{%
u}\left( .,t\right) \right\Vert _{L^{2}\left( O_{\frac{r}{8}};H\right) }+e^{%
\frac{-\gamma r^{2}}{64}}N_{\gamma }\text{ when }0\leq t\leq 1  \tag{6.10}
\end{equation}%
Then, $(6.8)-(6.10)$ show that there is a constant $N_{\gamma ,\varepsilon
,V}$, which such that 
\[
e^{C\left( \gamma ,\varepsilon \right) r^{2}}\left\Vert \tilde{u}\left(
.,0\right) \right\Vert _{X}\leq N_{\gamma ,\varepsilon ,V}. 
\]%
For $r\rightarrow \infty $ we obtain $u\equiv 0.$

\textbf{Proof of Theorem 4. }F\i rst all of, we show the following:

\textbf{Lemma 6.2. }Let the assumptions (1)-(2) of Conditon 3.1 hold.
Moreover, let%
\[
V\in B\text{, }\lim\limits_{r\rightarrow \infty }\left\Vert V\right\Vert
_{O\left( r\right) }=0. 
\]%
Then the estimate

\begin{equation}
r\sqrt{\frac{\varepsilon }{8\mu }}\left\Vert e^{\varkappa -\sigma +\chi
}\upsilon \right\Vert _{L^{2}\left( R^{n+1};H\right) }\leq \left\Vert
e^{\varkappa -\sigma +\chi }\left[ \partial _{t}u-\Delta u-Au\right]
\upsilon \right\Vert _{L^{2}\left( R^{n+1};H\right) }  \tag{6.11}
\end{equation}%
holds, when $\varepsilon >0$, $\mu >0$, $r>0$ and $\upsilon \in
C_{0}^{\infty }\left( R^{n+1};H\left( A\right) \right) $, where%
\[
\chi \left( t\right) =\frac{r^{2}t(1-t)\left( 1-2t\right) }{6}. 
\]

\textbf{Proof.} Let\textbf{\ }$f=e^{\varkappa +\chi -\sigma }\upsilon $.
Then, 
\[
e^{\varkappa +\chi -\sigma }\left[ \partial _{t}u-\left( \Delta u+Au\right) %
\right] \upsilon =\partial _{t}f-Sf-Kf. 
\]%
From $(3.8)-(3.10)$ with $\gamma =1$, $a+ib=1$ and $\varphi \left(
x,t\right) =\varkappa \left( x,t\right) +\chi \left( t\right) -\sigma \left(
\varepsilon ,t\right) $ we have 
\[
S=\Delta +A+4%
\mu
^{2}i|x+rt(1-t)e_{1}|^{2}+2%
\mu
ni+\text{ } 
\]%
\[
2%
\mu
R(1-2t)(x_{1}+rt(1-t))-\sigma +\left( t^{2}-t+\frac{1}{6}\right) r^{2}, 
\]%
\[
K=-4%
\mu
(x+rt(1-t)e_{1})\text{\textperiodcentered }\nabla -2%
\mu
n\text{, } 
\]%
\[
S_{t}+[S,K]=-8%
\mu
\triangle +32%
\mu
^{3}|x+rt(1-t)e_{1}|^{2}+ 
\]

\[
4%
\mu
r(4\mu \left( 1-2t-1\right) \left( (x_{1}+rt(1-t)\right) +\left( 2t-1\right)
r^{2}+\frac{\left( 1+\varepsilon \right) r^{2}}{8\mu } 
\]%
and

\begin{equation}
(S_{t}f+[S,K]f,f)_{X}=32%
\mu
^{3}\dint\limits_{R^{n}}\left\vert x+rt(1-t)e_{1}+\frac{(4\mu \left(
1-2t-1\right) r}{16\mu ^{2}}e_{1}\right\vert ^{2}\left\Vert f\right\Vert
^{2}dx+  \tag{6.12}
\end{equation}%
\[
8\mu \dint\limits_{R^{n}}\left\vert \nabla f\right\vert _{H}^{2}dx+\frac{%
\varepsilon r^{2}}{8\mu }\dint\limits_{R^{n}}\left\Vert f\right\Vert
^{2}dx\geq \frac{\varepsilon r^{2}}{8\mu }\dint\limits_{R^{n}}\left\Vert
f\right\Vert ^{2}dx. 
\]

Then from $\left( 6.12\right) $ a similar way as Lemma 6.1 we obtain the
estimate $\left( 6.11\right) .$

\bigskip \textbf{Proof. }Assume that $u$ verifies the conditions in Theorem
4 and let $\tilde{u}$ be the Appel transformation of $u$ defined in Lemma
4.1 with $a+ib=1$, $\alpha =1$ and $\beta =1+\frac{2}{\beta }$. $\tilde{u}%
\in L^{\infty }\left( 0,1;X\left( A\right) \right) \cap L^{2}\left(
0,1;Y^{1}\right) $ is a solution of the equation 
\[
\partial _{t}u=\Delta u+Au+\tilde{V}u,\text{ }x\in R^{n},\text{ }t\in \left[
0,1\right] 
\]%
with $\tilde{V}$ a bounded potential in $R^{n}\times \lbrack 0,1]$ and $%
\gamma =\frac{1}{2\delta }.$ Then, we have 
\[
\left\Vert e^{\gamma \left\vert x\right\vert ^{2}}\tilde{u}\left( .,0\right)
\right\Vert _{X}=\left\Vert \tilde{u}\left( .,0\right) \right\Vert _{X}\text{%
, }\left\Vert e^{\gamma \left\vert x\right\vert ^{2}}\tilde{u}\left(
.,1\right) \right\Vert _{X}=\left\Vert \tilde{u}\left( .,1\right)
\right\Vert _{X}. 
\]%
From Lemma 3.3 and Lemma 3. 4 with $a+ib=1$, we have

\[
\sup\limits_{t\in \left[ 0,1\right] }\left\Vert e^{\gamma \left\vert
x\right\vert ^{2}}\tilde{u}\left( .,t\right) \right\Vert _{X}+\left\Vert 
\sqrt{t\left( 1-t\right) }e^{\gamma \left\vert x\right\vert ^{2}}\nabla 
\tilde{u}\right\Vert _{L^{2}\left( R^{n}\times \left[ 0,1\right] ;H\right)
}\leq 
\]%
\[
e^{\left( M_{1}+M_{1}^{2}\right) }\left[ \left\Vert e^{\gamma \left\vert
x\right\vert ^{2}}\tilde{u}\left( .,0\right) \right\Vert _{X}+\left\Vert
e^{\gamma \left\vert x\right\vert ^{2}}\tilde{u}\left( .,1\right)
\right\Vert _{X}\right] , 
\]%
where 
\[
M_{1}=\left\Vert \tilde{V}\right\Vert _{B}. 
\]

The proof is finished by setting $\upsilon (x,t)=\theta _{M}(x)\eta _{R}(t)%
\tilde{u}(x,t)$, by using Carleman inequality $\left( 6.11\right) $ and in\
similar argument that we used to prove Theorem 1.

\begin{center}
\textbf{7. Unique continuation properties for the system of Schr\"{o}dinger
equations }
\end{center}

Consider the system of Schr\"{o}dinger equation%
\begin{equation}
\partial _{t}u_{m}=i\left[ \Delta
u_{m}+\sum\limits_{j=1}^{N}a_{mj}u_{j}+\sum\limits_{j=1}^{N}b_{mj}u_{j}%
\right] ,\text{ }x\in R^{n},\text{ }t\in \left( 0,T\right) ,  \tag{7.1}
\end{equation}%
where $u=\left( u_{1},u_{2},...,u_{N}\right) ,$ $u_{j}=u_{j}\left(
x,t\right) ,$ $a_{mj}$ \ are real-valued and $b_{mj}=b_{mj}\left( x,t\right) 
$ are complex valued functions$.$ Let $l_{2}=l_{2}\left( N\right) $ and $%
l_{2}^{s}=l_{2}^{s}\left( N\right) $ (see $\left[ \text{23, \S\ 1.18}\right] 
$). Let $A$ be the operator in $l_{2}\left( N\right) $ defined by%
\[
\text{ }D\left( A\right) =\left\{ u=\left\{ u_{j}\right\} ,\text{ }%
\left\Vert u\right\Vert _{D\left( A\right) }=\left(
\sum\limits_{m,j=1}^{N}\left( a_{mj}u_{j}\right) ^{2}\right) ^{\frac{1}{2}%
}<\infty \right\} , 
\]

\begin{equation}
A=\left[ a_{mj}\right] \text{, }m,j=1,2,...,N,\text{ }N\in \mathbb{N} 
\tag{7.2}
\end{equation}%
and

\[
V\left( x,t\right) =\left[ b_{mj}\left( x,t\right) \right] \text{, }%
m,j=1,2,...,N. 
\]

Let 
\[
X_{2}=L^{2}\left( R^{n};l_{2}\right) ,Y^{s,2}=H^{s,2}\left(
R^{n};l_{2}\right) . 
\]

From Theorem 1 we obtain the following result

\textbf{Theorem 7.1. }Assume:

(1) $a_{mj}=a_{jm},$ $a_{mj}\in C^{\left( 1\right) }\left( R^{n}\right) ,$ $%
\sum\limits_{m,j=1}^{N}a_{mj}>0,$ Moreover, 
\[
\dsum\limits_{k=1}^{n}\left( x_{k}\left[ A\frac{\partial f}{\partial x_{k}}-%
\frac{\partial A}{\partial x_{k}}f\right] ,f\right) _{X_{2}}\geq 0\text{ for 
}f\in C^{1}\left( R^{n};l_{2}^{s}\right) , 
\]%
$\func{Im}\sum\limits_{m,j=1}^{N}b_{mj}\left( x,t\right) >0$ for $x\in
R^{n}, $ $t\in \left[ 0,T\right] $ and\ $\func{Im}\sum%
\limits_{m,j=1}^{N}b_{mj}\left( x,t\right) \in L^{1}\left( 0,T;L^{\infty
}\left( R^{n}\right) \right) ;$

(2)$\sup\limits_{t\in \left[ 0,1\right] }\left\Vert e^{\left\vert
x\right\vert ^{2}\mu ^{-2}\left( t\right) }b_{mj}\left( .,t\right)
\right\Vert _{L^{\infty }\left( R^{n};B\left( l_{2}\right) \right) }<\infty
, $ where%
\[
\mu \left( t\right) =\alpha t+\beta \left( 1-t\right) \text{, }\alpha ,\beta
>0,\alpha \beta <2; 
\]

(3) $u\in C\left( \left[ 0,1\right] ;l_{2}\right) $ be a solution of the
equation $\left( 7.1\right) $ and%
\[
\left\Vert e^{\frac{\left\vert x\right\vert ^{2}}{\beta ^{2}}}u\left(
.,0\right) \right\Vert _{X_{2}}<\infty ,\left\Vert e^{\frac{\left\vert
x\right\vert ^{2}}{\alpha ^{2}}}u\left( .,T\right) \right\Vert
_{X_{2}}<\infty . 
\]

Then $u\left( x,t\right) \equiv 0.$

\ \textbf{Proof.} Consider the operators $A$ and $V\left( x,t\right) $ in $%
H=l^{2}$ defined by $\left( 7.2\right) $. Then the problem $\left(
8.1\right) -\left( 8.2\right) $ can be rewritten as the problem $\left(
1.1\right) $, where $u\left( x\right) =\left\{ u_{j}\left( x\right) \right\}
,$ $f\left( x\right) =\left\{ f_{j}\left( x\right) \right\} $, $j=1,2,...,N,$%
\ $x\in R^{n}$ are the functions with values in\ $H=l^{2}$. It is easy to
see that $A$ is a symmetric operator in $l_{2}$ and other conditions of
Theorem 1 are satisfied. Hence, from Teorem1 we obtain the conculision.

\begin{center}
\textbf{8. Unique continuation properties for anisotropic Schr\"{o}dinger
equation }

\ \ \ \ \ \ \ \ \ \ \ \ \ \ \ \ \ \ \ \ \ \ \ \ \ \ \ \ \ \ \ \ \ \ \ \ \ \
\ \ \ \ \ \ \ \ \ \ 
\end{center}

Let us consider the following problem

\begin{equation}
\partial _{t}u=i\left[ \Delta _{x}u+\sum\limits_{\left\vert \alpha
\right\vert =2m}a_{\alpha }D_{y}^{\alpha }u\left( x,y,t\right)
+\dint\limits_{G}K\left( x,y,\tau ,t\right) u\left( x,y,\tau ,t\right) d\tau %
\right] ,  \tag{8.1}
\end{equation}%
\[
\text{ }x\in R^{n},\text{ }y\in G,\text{ }t\in \left[ 0,T\right] ,\text{ } 
\]

\begin{equation}
B_{j}u=\sum\limits_{\left\vert \beta \right\vert \leq m_{j}}\ b_{j\beta
}D_{y}^{\beta }u\left( x,y,t\right) =0\text{, }x\in R^{n},\text{ }y\in
\partial G,\text{ }j=1,2,...,m,  \tag{8.2}
\end{equation}%
where $a_{\alpha }$ are real valued function on $G\subset $ $R^{d},$ $d\geq
2,$ $b_{j\beta }$ are the complex valued functions on $\partial G$, $\alpha
=\left( \alpha _{1},\alpha _{2},...,\alpha _{n}\right) $, $\beta =\left(
\beta _{1},\beta _{2},...,\beta _{n}\right) ,$ $\mu _{i}<2m,$ $K=K\left(
x,y,\tau ,t\right) $ is a complex valued bounded function in $\Omega \times
G\times \left[ 0,T\right] ,$ here $\Omega =R^{n}\times G$, $G$ is a bounded
domain with sufficiently smooth $\left( d-1\right) $-dimensional boundary $%
\partial G$ and 
\[
D_{x}^{k}=\frac{\partial ^{k}}{\partial x^{k}},\text{ }D_{j}=-i\frac{%
\partial }{\partial y_{j}},\text{ }D_{y}=\left( D_{1,}...,D_{d}\right) ,%
\text{ }y=\left( y_{1},...,y_{d}\right) , 
\]

$\bigskip $Let%
\[
\xi ^{\prime }=\left( \xi _{1},\xi _{2},...,\xi _{d-1}\right) \in R^{d-1},%
\text{ }\alpha ^{\prime }=\left( \alpha _{1},\alpha _{2},...,\alpha
_{d-1}\right) \in Z^{d-1},\text{ } 
\]%
\[
\text{ }A\left( \xi ^{\prime },D_{y}\right) =\sum\limits_{\left\vert \alpha
^{\prime }\right\vert +j\leq 2m}a_{\alpha ^{\prime }}\xi _{1}^{\alpha
_{1}}\xi _{2}^{\alpha _{2}}...\xi _{d-1}^{\alpha _{d-1}}D_{y}^{j}\text{ ,} 
\]%
\[
B_{j}\left( \xi ^{\prime },D_{y}\right) =\sum\limits_{\left\vert \beta
^{\prime }\right\vert +j\leq m_{j}}b_{j\beta ^{\prime }}\xi _{1}^{\beta
_{1}}\xi _{2}^{\beta _{2}}...\xi _{d-1}^{\beta _{d-1}}D_{y}^{j}. 
\]

\textbf{Theorem 8.1}. Let the following conditions be satisfied:

\bigskip (1) $\Omega \in C^{2}$, $a_{\alpha }\in C^{\left( 1\right) }\left( 
\bar{\Omega}\right) $ for $\left\vert \alpha \right\vert =2m,$ $a_{\alpha
}\in L^{\infty }\left( \Omega \right) $ for $\left\vert \alpha \right\vert
<2m,$ Moreover, 
\[
\dsum\limits_{k=1}^{n}\left( x_{k}\left[ A\frac{\partial f}{\partial x_{k}}-%
\frac{\partial A}{\partial x_{k}}f\right] ,f\right) _{L^{2}\left( \Omega
\right) }\geq 0\text{ for }f\in C^{1}\left( R^{n};W^{2m,2}\left( G\right)
\right) , 
\]%
$b_{j\beta }\in C\left( \partial \Omega \right) ,$ $\sum%
\limits_{j=1}^{m}b_{j\beta }\sigma _{j}\neq 0,$ for $\left\vert \beta
\right\vert =m_{j},$ where $\sigma =\left( \sigma _{1},\sigma
_{2},...,\sigma _{d}\right) \in R^{d}$ is a normal to $\partial G$ $;$

(2) for $\xi \in R^{n}$, $\lambda \in S\left( \varphi \right) =\left\{
\lambda \in \mathbb{C}\text{, }\left\vert \arg \lambda \right\vert \leq
\varphi ,\text{ }0\leq \varphi <\pi \right\} $ and $\left\vert \xi
\right\vert +\left\vert \lambda \right\vert \neq 0$ let $\lambda +$ $%
\sum\limits_{\left\vert \alpha \right\vert =2m}a_{\alpha }\xi ^{\alpha }\neq
0$;

(3) the problem%
\[
\lambda +A\left( \xi ^{\prime },D_{y}\right) \vartheta \left( y\right) =0,%
\text{ }B_{j}\left( \xi ^{\prime },D_{y}\right) \vartheta \left( 0\right)
=h_{j}\text{, }j=1,2,...,m 
\]

has a unique solution $\vartheta \in C_{0}\left( \mathbb{R}_{+}\right) $ for
all $h=\left( h_{1},h_{2},...,h_{d}\right) \in \mathbb{C}^{d}$ and for $\xi
^{\prime }\in R^{d-1};$

(4) $\func{Im}\dint\limits_{R^{n}}\dint\limits_{G}K\left( x,y,\tau ,t\right)
d\tau dy>0$ for $x\in R^{n}$ and $t\in \left[ 0,T\right] ,$ $T\in \left[ 0,1%
\right] ;$ moreover,\ 
\[
\func{Im}\dint\limits_{R^{n}}\dint\limits_{G}K\left( x,y,\tau ,t\right)
d\tau dy\in L^{1}\left( 0,T;L^{\infty }\left( G\right) \right) ; 
\]

(5) 
\[
\sup\limits_{t\in \left[ 0,1\right] }\left\Vert e^{\left\vert x\right\vert
^{2}\mu ^{2}\left( t\right) }K\left( .,.,t\right) \right\Vert _{L^{\infty
}\left( R^{n};L^{2}\left( G\right) \right) }<\infty ,\text{ } 
\]

where, $\mu \left( t\right) =\alpha t+\beta \left( 1-t\right) $, $\alpha
,\beta >0,\alpha \beta <2;$

(5) Assume $u\in C\left( \left[ 0,1\right] ;L^{p}\left( \Omega \right)
\right) $ be a solution of the equation $\left( 8.1\right) -\left(
8.2\right) $ and%
\[
\left\Vert e^{\frac{\left\vert x\right\vert ^{2}}{\beta ^{2}}}u\left(
.,.,0\right) \right\Vert _{L^{2}\left( \Omega \right) }<\infty ,\left\Vert
e^{\frac{\left\vert x\right\vert ^{2}}{\alpha ^{2}}}u\left( .,.,T\right)
\right\Vert _{L^{2}\left( \Omega \right) }<\infty . 
\]

Then $u\left( x,y,t\right) \equiv 0.$

\ \textbf{Proof. }Let us consider operators $A$ and $V\left( x,t\right) $ in 
$H=L^{2}\left( G\right) $ that are defined by the equalities 
\[
D\left( A\right) =\left\{ u\in W^{2m,2}\left( G\right) \text{, }B_{j}u=0,%
\text{ }j=1,2,...,m\text{ }\right\} ,\ Au=\sum\limits_{\left\vert \alpha
\right\vert =2m}a_{\alpha }D_{y}^{\alpha }u\left( y\right) , 
\]%
\[
V\left( x,t\right) u=\dint\limits_{G}K\left( x,y,\tau ,t\right) u\left(
x,y,\tau ,t\right) d\tau . 
\]

Then the problem $\left( 8.1\right) -\left( 8.2\right) $ can be rewritten as
the problem $\left( 1.1\right) $, where $u\left( x\right) =u\left(
x,.\right) ,$ $f\left( x\right) =f\left( x,.\right) $,\ $x\in \sigma $ are
the functions with values in\ $H=L^{2}\left( G\right) $. By virtue of $\left[
\text{1}\right] $ operator $A+\mu $ is positive in $L^{2}\left( G\right) $
for sufficiently large $\mu >0$. Moreover, in view of (1)-(5) all conditons
of Theorem 1 are hold. Then Theorem1 implies the assertion.

\begin{center}
\textbf{9.} \textbf{The Wentzell-Robin type mixed problem for Schr\"{o}%
dinger equations}
\end{center}

Consider the problem $\left( 1.5\right) -\left( 1.6\right) $. \ Let 
\[
\sigma =R^{n}\times \left( 0,1\right) . 
\]

In this section, we present the following result:

\bigskip \textbf{Theorem 9.1. } Suppose the the following conditions are
satisfied:

(1) let $a\left( x,.\right) $ be positive, $b\left( x,.\right) $ be a
real-valued function on $\left( 0,1\right) ,$ $a\in C^{\left( 1\right)
}\left( \left[ 0,1\right] \times R^{n}\right) $, Moreover, 
\[
\dsum\limits_{k=1}^{n}\left( x_{k}\left[ A\frac{\partial f}{\partial x_{k}}-%
\frac{\partial A}{\partial x_{k}}f\right] ,f\right) _{L^{2}\left( \sigma
\right) }\geq 0\text{ for }f\in C^{1}\left( R^{n};W^{2,2}\left( 0,1\right)
\right) , 
\]%
$b\left( .,y\right) \in C\left( R^{n}\right) $ for a.e. $y\in \left[ 0,1%
\right] ,$ $b\left( x,.\right) \in L_{\infty }\left( 0,1\right) $ for a.e. $%
x\in R^{n}$\ and%
\[
\exp \left( -\dint\limits_{\frac{1}{2}}^{x}b\left( \tau \right) a^{-1}\left(
x,\tau \right) d\tau \right) \in L_{1}\left( 0,1\right) \text{ for a.e.}x\in
R^{n}; 
\]

(2)\ $\func{Im}\dint\limits_{R^{n}}\dint\limits_{0}^{1}K\left( x,\tau
,t\right) d\tau dy>0$ for $x\in R^{n}$ and $t\in \left[ 0,T\right] $;
moreover,\ 
\[
\func{Im}\dint\limits_{R^{n}}\dint\limits_{0}^{1}K\left( x,\tau ,t\right)
d\tau dy\in L^{1}\left( 0,T;L^{\infty }\left( R^{n}\right) \right) ; 
\]

(3) 
\[
\sup\limits_{t\in \left[ 0,1\right] }\left\Vert e^{\left\vert x\right\vert
^{2}\mu ^{2}\left( t\right) }V\left( .,y,t\right) \right\Vert _{L^{\infty
}\left( \sigma \right) }<\infty \text{ for }y\in \left[ 0,1\right] , 
\]%
where, $\mu \left( t\right) =\alpha t+\beta \left( 1-t\right) $, $\alpha
,\beta >0,\alpha \beta <2;$

(4) $\ u\in C\left( \left[ 0,T\right] ;L^{2}\left( \sigma \right) \right) $
be a solution of the equation $\left( 1.5\right) -\left( 1.6\right) $ and%
\[
\left\Vert e^{\frac{\left\vert x\right\vert ^{2}}{\beta ^{2}}}u\left(
.,0\right) \right\Vert _{L^{2}\left( \sigma \right) }<\infty ,\left\Vert e^{%
\frac{\left\vert x\right\vert ^{2}}{\alpha ^{2}}}u\left( .,T\right)
\right\Vert _{L^{2}\left( \sigma \right) }<\infty . 
\]

Then $u\left( x,y,t\right) \equiv 0.$

\ \textbf{Proof.} Let us consider the operator $A$ in $H=L^{2}\left(
0,1\right) $ defined by $\left( 1.4\right) .$ Then $\left( 8.1\right)
-\left( 8.2\right) $ can be rewritten as the problem $\left( 1.1\right) $,
where $u\left( x\right) =u\left( x,.\right) ,$ $f\left( x\right) =f\left(
x,.\right) $,\ $x\in \sigma $ are the functions with values in\ $%
H=L^{2}\left( 0,1\right) $. By virtue of $\left[ \text{10, 11}\right] $ the
operator $A$ generates analytic semigroup in $L^{2}\left( 0,1\right) $.
Hence, by virtue of (1)-(3), all conditons of Theorem 1 are satisfied. Then
Theorem1 implies the assertion.

\textbf{References}\ \ 

\begin{enumerate}
\item H. Amann, Operator-valued Fourier multipliers, vector-valued Besov
spaces, and applications, Math. Nachr. 186 (1997), 5-56.\ 

\item H. Amann, Linear and quasi-linear equations,1, Birkhauser, Basel 1995.

\item A. Bonami, B. Demange, A survey on uncertainty principles related to
quadratic forms. Collect. Math. 2006, Vol. Extra, 1--36.

\item C. E. Kenig, G. Ponce, L. Vega, On unique continuation for nonlinear
Schr\"{o}dinger equations, Comm. Pure Appl. Math. 60 (2002) 1247--1262.

\item L. Escauriaza, C. E. Kenig, G. Ponce, L. Vega, On uniqueness
properties of solutions of Schr\"{o}dinger Equations, Comm. PDE. 31, 12
(2006) 1811--1823.

\item L. Escauriaza, C. E. Kenig, G. Ponce, L. Vega, On uniqueness
properties of solutions of the k-generalized KdV, J. of Funct. Anal. 244, 2
(2007) 504--535.

\item L. Escauriaza, C. E. Kenig, G. Ponce, and L. Vega, Hardy's uncertainty
principle, convexity and Schr\"{o}dinger Evolutions, J. European Math. Soc.
10 (4) (2008) 883--907.

\item L. H\"{o}rmander, A uniqueness theorem of Beurling for Fourier
transform pairs, Ark. Mat. 29, 2 (1991) 237--240.

\item J. A. Goldstain, Semigroups of Linear Operators and Applications,
Oxford University Press, Oxfard, 1985.

\item A. Favini, G. R. Goldstein, J. A. Goldstein and S. Romanelli,
Degenerate Second Order Differential Operators Generating Analytic
Semigroups in $L_{p}$ and $W^{1,p}$, Math. Nachr. 238 (2002), 78 --102.

\item V. Keyantuo, M. Warma, The wave equation with Wentzell--Robin boundary
conditions on Lp-spaces, J. Differential Equations 229 (2006) 680--697.

\item P. Guidotti, \ Optimal regularity for a class of singular abstract
parabolic equations, J. Differential Equations, v. 232, 2007, 468--486.

\item Lunardi A., Analytic Semigroups and Optimal Regularity in Parabolic
Problems, Birkhauser, 2003.

\item Lions, J-L., Magenes, E., Nonhomogenous Boundary Value Broblems, Mir,
Moscow, 1971.

\item N. Okazawa, T. Suzuki, T. Yokota, Energy methods for abstract
nonlinear Schredinger equations, Evol. Equ. Control Theory 1(2)(2012),
337-354.

\item V. B. Shakhmurov, Nonlinear abstract boundary value problems in
vector-valued function spaces and applications, Nonlinear Anal-Theor., v.
67(3) 2006, 745-762.

\item V. B. Shakhmurov, Coercive boundary value problems for regular
degenerate differential-operator equations, J. Math. Anal. Appl., 292 ( 2),
(2004), 605-620.

\item R. Shahmurov, On strong solutions of a Robin problem modeling heat
conduction in materials with corroded boundary, Nonlinear Anal., Real World
Appl., v.13, (1), 2011, 441-451.

\item R. Shahmurov, Solution of the Dirichlet and Neumann problems for a
modified Helmholtz equation in Besov spaces on an annuals, J. Differential
equations, v. 249(3), 2010, 526-550.

\item \ E. M. Stein, R. Shakarchi, Princeton, Lecture in Analysis II.
Complex Analysis, Princeton, University Press (2003).

\item A. Sitaram, M. Sundari, S. Thangavelu, Uncertainty principles on
certain Lie groups, Proc. Indian Acad. Sci. Math. Sci. 105 (1995), 135-151.

\item L. Weis, Operator-valued Fourier multiplier theorems and maximal $%
L_{p} $ regularity, Math. Ann. 319, (2001), 735-758.

\item H. Triebel, Interpolation theory, Function spaces, Differential
operators, North-Holland, Amsterdam, 1978.

\item S. Yakubov and Ya. Yakubov, \textquotedblright Differential-operator
Equations. Ordinary and Partial \ Differential Equations \textquotedblright
, Chapman and Hall /CRC, Boca Raton, 2000.
\end{enumerate}

\end{document}